\documentstyle[amscd, amstex]{amsart}

\def\ZZ         {{\Bbb Z}}
\def\RR         {{\Bbb R}}
\def\CC         {{\Bbb C}}
\def\QQ         {{\Bbb Q}}
\def\PP         {{\Bbb P}}
\def\TT         {{\Bbb T}}
\def\ii         {{\rm i}}
\def\ee         {{\rm e}}
\def\Ext        {{\rm Ext}}
\def\Hom        {{\rm Hom}}

\def\rk         {{\rm rk}}
\def\st         {{\rm st}}

\def\iso        {{\cong}}
\def\contr      {{\lrcorner}}

\def\qed        {{\hfill{$\Box$}}}

\def\dim        {{\rm dim}}

\def\deg        {{\rm deg}}
\def\codim      {{\rm codim}}

\def\Box        {{\square }}
\def\cal        {\mathcal}

\newtheorem{thm}{Theorem}[section]
\newtheorem{lem}[thm]{Lemma}
\newtheorem{cor}[thm]{Corollary}
\newtheorem{pr}[thm]{Proposition}

\theoremstyle{definition}
\newtheorem{rem}[thm]{Remark}

\newtheorem{defn}[thm]{Definition}
\newtheorem{conj}[thm]{Conjecture}
\newtheorem{cdefn}[thm]{Conjectural Definition}

\newcommand{\psx}{{\PP}_{\Sigma_X}}

\newcommand{\ts}{{\TT}_\sigma}

\newcommand{\key}{\bibitem}

\newcommand{\ps}{{{\PP}_{\Sigma}}}

\newcommand{\gr}{{\rm Gr}}

\newcommand\hidot{{\raise1pt\hbox{$\scriptscriptstyle\bullet$}}}
\newcommand\lodot{{\raise.3pt\hbox{$\scriptscriptstyle\bullet$}}}

\newcommand{\pp}{{\PP}}

\newcommand{\inte}{{\rm int}}

\begin{document}

\title
[String cohomology of Calabi-Yau hypersurfaces]
{String cohomology of Calabi-Yau hypersurfaces via Mirror Symmetry}

\author{Lev A. Borisov}
\address{Department of Mathematics, Columbia University,
New York, NY  10027, USA}
\email{lborisov@@math.columbia.edu}

\author{Anvar R. Mavlyutov}
\address {Max-Planck-Institut f\"ur Mathematik, Bonn, D-53111, Germany}
\email{anvar@@mpim-bonn.mpg.de}


\begin{abstract}
We propose a construction of string cohomology spaces for
Calabi-Yau hypersurfaces that arise in Batyrev's mirror symmetry
construction. The spaces are defined  explicitly in terms of the
corresponding reflexive polyhedra in a mirror-symmetric manner. We
draw connections with other approaches to the string cohomology,
in particular with the work of Chen and Ruan.
\end{abstract}

\maketitle

\tableofcontents

\section{Introduction}\label{section.intro}

The notion of orbifold cohomology has appeared in physics
as a result of studying the string theory on orbifold global quotients,
(see \cite{dhvw}). In addition to the usual cohomology of the quotient,
this space was supposed to include the so-called twisted sectors,
whose existence was predicted by the modular invariance condition
on the partition function of the theory.
Since then,
there have been several attempts to give a rigorous mathematical formulation of this
cohomology theory. The first two, due to \cite{bd} and \cite{Batyrev.cangor},
 tried to define the topological
invariants of certain algebraic varieties (including orbifold
global quotients) that should correspond to the dimensions of the
Hodge components of a conjectural string cohomology space. These
invariants should have the property arising naturally from
physics: they are preserved by partial crepant resolutions;
moreover, they coincide with the usual Hodge numbers for smooth
varieties. Also, these invariants must be the same as those
defined by physicists for orbifold global quotients. In
\cite{Batyrev.cangor, Batyrev.nai}, Batyrev has successfully
solved this problem for a large class of singular algebraic
varieties. The first mathematical definition of the orbifold
cohomology {\em space} was given in \cite{cr} for arbitrary
orbifolds. Moreover, this orbifold cohomology possesses a  product
structure arising as a limit of a natural quantum product. It is
still not entirely clear if the dimensions of the Chen-Ruan
cohomology coincide with the prescription of Batyrev whenever both
are defined, but they do give the same result for reduced global
orbifolds.

In this paper, we propose a construction of string cohomology
spaces for Calabi-Yau hypersurfaces that arise in the Batyrev
mirror symmetry construction (see \cite{b2}), with the spaces
defined rather explicitly in terms of the corresponding reflexive
polyhedra. A peculiar feature of our construction is that instead
of a single string cohomology space we construct a
finite-dimensional family of such spaces, which is consistent with
the physicists' picture (see \cite{Greene}). We verify that this
construction is consistent with the previous definitions in
\cite{bd}, \cite{Batyrev.cangor} and \cite{cr}, in the following
sense. The (bigraded) dimension of our space coincides with the
definitions of \cite{bd} and \cite{Batyrev.cangor}. In the case of
hypersurfaces that have only orbifold singularities, we recover
Chen-Ruan's orbifold cohomology as one special element of this
family of string cohomology spaces. We also conjecture a partial
natural ring structure on our string cohomology space, which is in
correspondence with the cohomology ring of crepant resolutions.
This may be used as a real test of the Chen-Ruan orbifold
cohomology ring. We go further, and conjecture the B-model chiral
ring on the string cohomology space. This is again consistent with
the description of the B-model chiral ring of smooth Calabi-Yau
hypersurfaces in \cite{m2}.

Our construction of the string cohomology space for Calabi-Yau hypersurfaces is
motivated by Mirror Symmetry. Namely, the description in \cite{m3} of
the cohomology of semiample hypersurfaces in toric varieties applies to the smooth Calabi-Yau
hypersurfaces in \cite{b2}. Analysis of Mirror Symmetry on this cohomology leads to a natural
construction of the string cohomology space for all semiample Calabi-Yau hypersurfaces.
As already mentioned, our string cohomology space depends not only
on the complex  structure (the  defining polynomial $f$), but also on some extra parameter
we call $\omega$. For special values of this parameter of an orbifold Calabi-Yau hypersurface,
we get the orbifold Dolbeault cohomology of \cite{cr}. However, for non-orbifold
Calabi-Yau hypersurfaces, there is no natural special choice of $\omega$,
which means that the general definition of the string cohomology space should depend on some
mysterious extra parameter. In the situation of Calabi-Yau hypersurfaces, the parameter $\omega$
corresponds to the defining polynomial of the mirror Calabi-Yau hypersurface.
In general, we expect that this parameter should be related to the ``stringy complexified
K\"ahler class'', which is yet to be defined.

In an attempt to extend our definitions beyond the Calabi-Yau hypersurface case,
we give a conjectural definition of string cohomology vector spaces
for stratified varieties with $\QQ$-Gorenstein toroidal singularities that satisfy certain
restrictions on the types of singular strata. This definition involves intersection
cohomology of the closures of strata, and we check that it produces spaces of correct bigraded
dimension. It also reproduces orbifold cohomology of a $\QQ$-Gorenstein
toric variety as a special case.

Here is an outline of our paper. In Section \ref{s:sth}, we examine the connection between
the original definition of the {\em string-theoretic} Hodge numbers in \cite{bd}
and the {\em stringy} Hodge numbers in \cite{Batyrev.cangor}. We point out that these do not
always give the same result and argue that the latter definition is the more useful one.
In Section \ref{section.mirr}, we briefly review the mirror symmetry construction
of Batyrev, mainly to fix our notations and to describe the properties we will use in
the derivation of the string cohomology. Section \ref{section.anvar} describes the cohomology
of semiample hypersurfaces in toric varieties and explains how mirror symmetry provides
a conjectural definition of the string cohomology of Calabi-Yau hypersurfaces.
It culminates in  Conjecture \ref{semiampleconj}, where we define the stringy Hodge spaces
of semiample Calabi-Yau hypersurfaces in complete toric varieties.
We spend most of the remainder of the paper establishing the expected properties of
the string cohomology space. Sections \ref{s:hd} and \ref{section.brel} calculate the
dimensions of the building blocks of our cohomology spaces. In Section \ref{section.brel},
we develop a theory of deformed semigroup rings which may be of independent interest.
This allows us to show in Section \ref{section.bbo}
that Conjecture \ref{semiampleconj} is compatible with the
definition of the stringy Hodge numbers from
\cite{Batyrev.cangor}. In the non-simplicial case, this requires the use
of $G$-polynomials of Eulerian posets, whose relevant properties are collected in the
Appendix.
Having established that the dimension is correct, we try to extend our
construction to the non-hypersurface case. Section \ref{section.general} gives another
conjectural definition of the string cohomology vector space in a somewhat more general
situation. It hints that the intersection cohomology and the perverse sheaves should play a
prominent role in future definitions of string cohomology. In Section
\ref{s:vs}, we connect our work with that of Chen-Ruan \cite{cr} and Poddar \cite{p}.
Finally, in Section \ref{section.vertex}, we provide yet another
description of the string cohomology of Calabi-Yau hypersurfaces, which
was inspired by the vertex algebra approach to Mirror Symmetry.

{\it Acknowledgments.} We thank Victor Batyrev, Robert Friedman,
Mainak Poddar, Yongbin Ruan and Richard Stanley for helpful
conversations and useful references. The second author also thanks
the Max-Planck Institut f\"ur Mathematik in Bonn for its
hospitality and support.

\section{String-theoretic and stringy Hodge numbers}
\label{s:sth}

The {\em string-theoretic} Hodge numbers were first defined in the
paper of Batyrev and Dais (see \cite{bd}) for varieties with
Gorenstein toroidal or quotient singularities. In subsequent
papers \cite{Batyrev.cangor,Batyrev.nai} Batyrev defined {\em
stringy} Hodge numbers for arbitrary varieties with log-terminal
singularities. To our knowledge, the relationship between these
two concepts has never been clarified in the literature. The goal
of this section is to show that the string-theoretic Hodge numbers
coincide with the stringy ones under some conditions on the
singular strata.

We begin with the definition of the string-theoretic Hodge numbers.
\begin{defn}\label{d:bd} \cite{bd}
Let $X = \bigcup_{i \in I} X_i$ be a stratified algebraic variety
over $\CC$ with at most Gorenstein toroidal singularities such
that for each $i \in I$ the singularities of $X$ along the stratum
$X_i$ of codimension $k_i$ are defined by a $k_i$-dimensional
finite rational polyhedral cone $\sigma_i$; i.e.,  $X$ is locally
isomorphic to $${\CC}^{k-k_i} \times U_{\sigma_i}$$ at each  point
$x \in X_i$,  where $U_{\sigma_i}$ is a $k_i$-dimensional  affine
toric variety which is associated with the cone $\sigma_i$ (see
\cite{d}), and $k=\dim X$. Then the polynomial $$ E^{\rm BD}_{\rm
st}(X;u,v) := \sum_{i \in I} E(X_i;u,v) \cdot S(\sigma_i,uv) $$ is
called the {\em string-theoretic E-polynomial} of $X$. Here,
$$S(\sigma_i,t):=(1-t)^{\dim \sigma_i}\sum_{n\in \sigma_i} t^{\deg
n}= (t-1)^{\dim \sigma_i}\sum_{n\in {\rm int} \sigma_i} t^{-\deg
n}$$ where $\deg$ is the linear function on $\sigma_i$ that takes
value $1$ on the generators of one-dimensional faces of
$\sigma_i$, and ${\rm int}\sigma_i$ is the relative interior of
$\sigma_i$. If we write $E_{\rm st}(X; u,v)$ in the form $$ E^{\rm
BD}_{\rm st}(X;u,v) = \sum_{p,q} a_{p,q} u^{p}v^{q}, $$ then the
numbers $h^{p,q{\rm (BD)}}_{\rm st}(X) := (-1)^{p+q}a_{p,q}$ are
called the {\em string-theoretic Hodge numbers} of $X$.
\end{defn}

\begin{rem}\label{r:epol}
 The E-polynomial in the above definition is defined for an arbitrary algebraic
variety $X$ as
$$E(X;u,v)=\sum_{p,q}e^{p,q}u^p v^q,$$
where $e^{p,q}=\sum_{k\ge0}(-1)^k h^{p,q}(H^k_c(X))$.
\end{rem}

Stringy Hodge numbers of $X$ are defined in terms of the resolutions of
its singularities. In general, one can only define the $E$-function in this
case, which may or not be a polynomial.
We refer to \cite{Kollar} for the definitions
of log-terminal singularities and related issues.

\begin{defn}\label{d:bcangor}\cite{Batyrev.cangor}
Let $X$ be a normal irreducible algebraic variety with at worst
log-terminal singularities, $\rho\, : \, Y \rightarrow
X$ a resolution of singularities such that the irreducible components
$D_1, \ldots, D_r$ of the exceptional locus is a divisor with simple
normal crossings. Let $\alpha_j>-1$ be the discrepancy of $D_j$,
see \cite{Kollar}.
Set $I: = \{1, \ldots, r\}$. For any
subset $J \subset I$ we
consider
\[ D_J := \left\{ \begin{array}{ll}
\bigcap_{ j \in J} D_j & \mbox{\rm if $J \neq \emptyset$}
\\
Y  & \mbox{\rm if $J =  \emptyset$} \end{array} \right. \;\;\;\;
\;\;\;\; \,\mbox{\rm and}   \;\;\;\; \;\;\;\; D_J^{\circ} := D_J
\setminus \bigcup_{ j \in\, I \setminus J} D_j. \] We define an
algebraic  function $E_{\rm st}(X; u,v)$ in two variables $u$ and
$v$ as follows: $$ E_{\rm st}(X; u,v) := \sum_{J \subset I}
E(D_J^{\circ}; u,v) \prod_{j \in J} \frac{uv-1}{(uv)^{a_j +1} -1}
$$ (it is  assumed   $\prod_{j \in J}$ to be  $1$, if $J =
\emptyset$). We call  $E_{\rm st}(X; u,v)$ {\em the stringy
$E$-function of} $X$. If $E_{\rm st}(X; u,v)$ is a polynomial,
define the stringy Hodge numbers the same way as
Definition~\ref{d:bd} does.
\end{defn}

It is not obvious at all that the above definition is independent
of the choice of the resolution. The original proof of Batyrev
uses a motivic integration over the spaces of arcs to relate the
$E$-functions obtained via different resolutions. Since the work
of D.~Abramovich, K. Karu, K. Matsuki, J. W\l odarsczyk
\cite{AKMW}, it is now possible to check the independence from the
resolution by looking at the case of a single blowup with a smooth
center compatible with the normal crossing condition.

\begin{lem}\label{strataE}
Let $X$ be a disjoint union of strata $X_i$, which are locally closed in
Zariski topology, and let $\rho$ be a resolution as in Definition~\ref{d:bcangor}.
For each $X_i$ consider
$$ E_{\rm st}(X_i\subseteq X; u,v) := \sum_{J \subset I}
E(D_J^{\circ}\cap \rho^{-1}(X_i); u,v)
\prod_{j \in J}
\frac{uv-1}{(uv)^{a_j +1} -1}. $$
Then this $E$-function is independent of the choice of the resolution
$Y$.
The $E$-function of $X$ decomposes as
$$
E_{\rm st}(X;u,v) = \sum_i  E_{\rm st}(X_i\subseteq X; u,v).
$$
\end{lem}

{\em Proof.}
Each resolution of $X$ induces a resolution of the complement of
$\bar{X_i}$. This shows that for each $X_i$
the sum
$$
\sum_{j,X_j\subseteq \bar{X_i}} E_{\rm st}(X_j\subseteq X; u,v)
$$
is independent from the choice of the resolution and is thus well-defined.
Then one uses the induction on dimension of $X_i$.
The last statement is clear.
\qed

\begin{rem}
It is a delicate question what data are really necessary to calculate
$E_{\rm st}(X_i\subseteq X; u,v)$. It is clear that the knowledge of
a Zariski open set of $X$ containing $X_i$ is enough. However, it
is not clear whether it is enough to know an analytic neighborhood
of $X_i$.
\end{rem}

We will use the above lemma to show that the string-theoretic Hodge numbers
and the stringy Hodge numbers coincide in a wide class of examples.

\begin{pr}\label{BDvsB}
Let $X=\bigcup_i X_i$ be a stratified algebraic variety with at worst
Gorenstein toroidal singularities as in Definition \ref{d:bd}.
Assume in addition that for each $i$ there is a desingularization $Y$ of
$X$ so that its restriction to the preimage of $X_i$
is a locally trivial fibration in Zariski topology. Moreover, for
a point $x\in X_i$ the preimage in $Y$ of an analytic neighborhood of $x$ is
complex-analytically isomorphic to a preimage of a neighborhood of $\{0\}$
in $U_{\sigma_i}$ under some resolution of singularities of $U_{\sigma_i}$,
times a complex disc, so that the isomorphism is compatible with the resolution
morphisms.
Then
$$
E_{\rm st}^{\rm BD}(X;u,v) = E_{\rm st}(X;u,v).
$$
\end{pr}

{\em Proof.}
Since $E$-polynomials are multiplicative for Zariski locally trivial
fibrations (see \cite{dk}),
the above assumptions on the singularities show that
$$
E_{\rm st}(X_i\subseteq X; u,v) =
E(X_i;u,v)E_{\rm st}(\{0\}\subseteq U_{\sigma_i};u,v).
$$
We have also used here the fact that since the fibers are projective,
the analytic isomorphism implies the algebraic one, by GAGA.
By the second statement of Lemma \ref{strataE}, it is enough to show that
$$
E_{\rm st}(\{0\}\subseteq U_{\sigma_i};u,v) = S(\sigma_i,uv).
$$
This result follows from the proof of \cite{Batyrev.cangor}, Theorem 4.3
where the products
$$
\prod_{j \in J}
\frac{uv-1}{(uv)^{a_j +1} -1}
$$
are interpreted as a geometric series and then as
sums of $t^{\deg(n)}$ over points $n$ of $\sigma_i$.
\qed

\begin{cor}
String-theoretic and stringy Hodge numbers coincide for
nondegenerate hypersurfaces (complete intersections) in Gorenstein
toric varieties.
\end{cor}

\begin{pf}
Indeed, in this case, the toric desingularizations of the ambient
toric variety induce the desingularizations with the required
properties.
\end{pf}

We will keep this corollary in mind and from now on will silently transfer
all the results on string-theoretic Hodge numbers of hypersurfaces and
complete intersections in toric varieties in \cite{bb}, \cite{bd}
to their stringy counterparts.

\begin{rem}
An example of the variety where string-theoretic and stringy Hodge numbers
{\em differ} is provided by the quotient of $\CC^2\times E$ by the finite group
of order six generated by
$$r_1:(x,y;z)\mapsto (x\ee^{2\pi\ii/3},y\ee^{-2\pi\ii/3};z),~
r_2:(x,y;z)\mapsto(y,x;z+p)$$
where $(x,y)$ are coordinates on $\CC^2$, $z$ is the uniformizing coordinate
on the elliptic curve $E$ and $p$ is a point of order two on $E$.
In its natural stratification, the quotient has a stratum of $A_2$
singularities, so that going around a loop in the stratum results in
the non-trivial automorphism of the singularity.
\end{rem}

\begin{rem}
We expect that the stringy Hodge numbers of algebraic varieties
with abelian quotient singularities coincide with the dimensions
of their orbifold cohomology, \cite{cr}. This is not going to be
true for the string-theoretic Hodge numbers. Also, the latter
numbers are not preserved by the partial crepant resolutions as
required by physics, see the above example. As a result, we
believe that the stringy Hodge numbers are the truly interesting
invariant, and that the string-theoretic numbers is a now obsolete
first attempt to define them.
\end{rem}

\section{Mirror symmetry construction of Batyrev}\label{section.mirr}

In this section, we review the mirror symmetry construction from \cite{b2}.
We can describe it starting with a semiample nondegenerate (transversal to the
torus orbits) anticanonical  hypersurface $X$ in a complete simplicial toric variety $\ps$.
Such a hypersurface is Calabi-Yau.
The semiampleness property produces a contraction map, the unique properties
of which are characterized by the following statement.

\begin{pr}\label{p:sem} \cite{m1}
Let $\PP_\Sigma$ be a complete toric variety with a big and nef
 divisor class $[X]\in A_{d-1}(\ps)$.
Then,
there exists a unique complete toric variety $\psx$ with a toric birational map
$\pi:\ps@>>>\psx$, such that   $\Sigma$ is a subdivision of
$\Sigma_X$, $\pi_*[X]$ is ample and $\pi^*\pi_*[X]=[X]$. Moreover, if
$X=\sum_{\rho}a_\rho D_\rho$ is torus-invariant, then $\Sigma_X$ is the normal fan of the associated polytope
$$\Delta_X=\{m\in M:\langle m,e_\rho\rangle\geq-a_\rho
\text{ for all } \rho\}\subset M_{\Bbb R}.$$
\end{pr}

\begin{rem} Our notation is a standard one taken from \cite{bc,c2}:
$M$ is a lattice of rank $d$;
$N=\text{Hom}(M,{\Bbb Z})$ is the dual lattice;
$M_{\Bbb R}$ and $N_{\Bbb R}$ are the $\Bbb R$-scalar extensions of $M$ and $N$;
$\Sigma$ is a finite rational polyhedral   fan in $N_{\Bbb R}$;
${\PP}_{\Sigma}$ is a  $d$-dimensional toric variety associated with   $\Sigma$;
$\Sigma(k)$ is the set of all $k$-dimensional cones in $\Sigma$;
$e_\rho$ is the minimal integral generator of
the $1$-dimensional
cone  $\rho\in\Sigma$ corresponding to a torus invariant irreducible divisor $D_\rho$.
\end{rem}

Applying Proposition~\ref{p:sem} to the semiample Calabi-Yau hypersurface,
we get that
the push-forward $\pi_*[X]$ is  anticanonical and ample,
whence, by Lemma 3.5.2
in \cite{ck}, the toric variety $\psx$ is
Fano, associated with the polytope
$\Delta\subset M_{\Bbb R}$
of the anticanonical divisor $\sum_{\rho} D_{\rho}$ on $\ps$.
Then, \cite[Proposition~2.4]{m1}
  shows that the image $Y:=\pi(X)$ is an ample nondegenerate hypersurface
in $\psx=\pp_\Delta$.
The fact that  $\pp_\Delta$ is Fano means by Proposition 3.5.5 in \cite{ck} that the polytope
$\Delta$ is reflexive,
i.e., its dual
$$\Delta^*=\{n\in N_{\Bbb R}: \langle m,n\rangle\ge-1\text{ for } m\in\Delta\}$$
has all its vertices at lattice points in $N$, and the only  lattice point in the interior of $\Delta^*$ is the origin $0$.
Now, consider  the  toric variety
$\pp_{\Delta^*}$ associated to the polytope $\Delta^*$ (the minimal
integral generators of its fan are precisely the vertices of $\Delta$).
 Theorem 4.1.9 in \cite{b2} says that
an anticanonical nondegenerate hypersurface
$Y^*\subset\pp_{\Delta^*}$ is a Calabi-Yau variety with canonical singularities.
The Calabi-Yau hypersurface $Y^*$ is expected to be a mirror of $Y$.
In particular,
they pass the topological mirror symmetry test for the stringy Hodge numbers:
$$h^{p,q}_{\rm st}(Y)=h_{\rm st}^{d-1-p,q}(Y^*), 0\le p,q\le d-1,$$
by \cite[Theorem~4.15]{bb}.
Moreover,  all crepant partial resolutions $X$  of   $Y$ have the same
stringy Hodge numbers:
$$h^{p,q}_{\rm st}(X)=h^{p,q}_{\rm st}(Y).$$
Physicists predict that  such resolutions of Calabi-Yau varieties
 have indistinguishable physical  theories.
Hence, all crepant partial  resolutions of $Y$ may be called the
 mirrors of crepant partial
resolutions of $Y^*$.
To connect this to the classical formulation of  mirror symmetry, one needs to note
that if there exist
  crepant smooth  resolutions $X$ and $X^*$ of $Y$ and $Y^*$, respectively,
then
$$h^{p,q}(X)=h^{d-1-p,q}(X^*), 0\le p,q\le d-1,$$
since the  stringy Hodge numbers coincide with the usual ones for smooth
Calabi-Yau varieties. The equality of Hodge numbers is expected
to extend to an isomorphism ({\it mirror map})
of the corresponding Hodge spaces, which is compatible with the chiral ring
products of A and B models (see \cite{ck} for more details).

\section{String cohomology construction for  Calabi-Yau hypersurfaces}
\label{section.anvar}

In this section, we show how the description of cohomology of semiample
hypersurfaces in \cite{m3} leads to a  construction
of the string cohomology space of  Calabi-Yau hypersurfaces.
We first review the building blocks participating in the description
of the cohomology in \cite{m3},
and then explain how these building blocks should interchange under mirror
symmetry for a pair of smooth Calabi-Yau hypersurfaces in Batyrev's mirror
symmetry construction. Mirror symmetry and the fact that the dimension of the
string cohomology
is the same for all partial crepant resolutions of ample
Calabi-Yau hypersurfaces leads us to a conjectural description
of string cohomology for all semiample Calabi-Yau hypersurfaces.
In the next three sections,
we will  prove that this space has the  dimension prescribed
by \cite{bd}.

The cohomology of a semiample nondegenerate hypersurface $X$ in a complete
simplicial toric variety $\ps$ splits into the {\it toric} and {\it residue} parts:
$$H^*(X)=H^*_{\rm toric}(X)\oplus H^*_{\rm res}(X),$$
where the first part is the image of the cohomology of the ambient space,
while the second is the residue map
image of the cohomology of the complement to the hypersurface.
By \cite[Theorem~5.1]{m2},
\begin{equation}\label{e:ann}
H^*_{\rm toric}(X)\cong H^*(\ps)/Ann(X)
\end{equation}
where $Ann(X)$ is the annihilator of the class $[X]\in H^2(\ps)$.
The cohomology of $\ps$  is isomorphic to
$${\Bbb C}[D_\rho:\rho\in\Sigma(1)]/(P(\Sigma)+SR(\Sigma)),$$
 where
$$P(\Sigma)=\biggl\langle \sum_{\rho\in\Sigma(1)}\langle m,e_\rho\rangle D_\rho:
m\in M\biggr\rangle$$
is the ideal of linear relations among the divisors,
and
$$SR(\Sigma)=\bigl\langle D_{\rho_1}\cdots D_{\rho_k}:\{e_{\rho_1},\dots,e_{\rho_k}\}
\not\subset\sigma
\text{ for all }\sigma\in\Sigma\bigr\rangle$$
is the Stanley-Reisner ideal.
Hence, $H^{*}_{\rm toric}(X)$ is isomorphic to  the bigraded ring
$$T(X)_{*,*}:={\Bbb C}[D_\rho:\rho\in\Sigma(1)]/I,$$
where $I=(P(\Sigma)+SR(\Sigma)):[X]$ is the ideal quotient,
and $D_\rho$ have the degree $(1,1)$.

The following modules over the ring $T(X)$ have appeared in the description
of cohomology of semiample hypersurfaces:

\begin{defn}
Given a big and nef class $[X]\in A_{d-1}(\ps)$   and $\sigma\in\Sigma_X$,
let
$$U^\sigma(X)=\biggl\langle \prod_{\rho\subset\gamma\in\Sigma}D_\rho:
\inte\gamma\subset\inte\sigma\biggr\rangle$$
be the bigraded ideal in ${\Bbb C}[D_\rho:\rho\in\Sigma(1)]$,
where $D_\rho$ have the degree (1,1).
Define the bigraded space
$$T^\sigma(X)_{*,*}=U^\sigma(X)_{*,*}/I^\sigma,$$
where
$$I^\sigma=\{u\in U^\sigma(X)_{*,*}:\,uvX^{d-\dim\sigma}
\in(P(\Sigma)+SR(\Sigma))\text{ for }v\in U^\sigma(X)_{\dim\sigma-*,\dim\sigma-*}\}.$$
\end{defn}

Next, recall from \cite{c} that
{\it any} toric variety $\ps$ has a homogeneous coordinate ring
$$S(\ps)={\Bbb C}[x_\rho:\rho\in\Sigma(1)]$$
 with variables $x_\rho$
corresponding to the irreducible torus invariant divisors   $D_\rho$.
This ring is graded by the Chow group $A_{d-1}(\ps)$, assigning $[\sum_{\rho} a_\rho D_\rho]$
to $\deg(\prod_{\rho} x_\rho^{a_\rho})$.
For a Weil divisor $D$ on $\ps$, there is an isomorphism
$H^0(\ps, O_\ps(D))\cong S(\ps)_\alpha$, where
$\alpha=[D]\in A_{d-1}(\ps)$. If $D$ is torus invariant, the monomials in
$S(\ps)_\alpha$ correspond to
the lattice points of the associated polyhedron $\Delta_D$.

In \cite{bc}, the following rings have been used to describe
the residue part of
cohomology of ample hypersurfaces in complete simplicial toric varieties:

\begin{defn}\label{d:r1} \cite{bc}
Given  $f\in S(\ps)_\beta$, set
$J_0(f):=\langle x_\rho\partial f/\partial x_\rho:\rho\in\Sigma(1)\rangle$ and
$J_1(f):=J_0(f):x_1\cdots x_n$.
Then define the rings
$R_0(f)=S(\ps)/J_0(f)$ and $R_1(f)=S(\ps)/J_1(f)$, which are
graded by the Chow group $A_{d-1}(\ps)$.
\end{defn}

In \cite[Definition~6.5]{m3}, similar rings were introduced to describe
the residue part of cohomology of semiample hypersurfaces:

\begin{defn}\label{d:rs1} \cite{m3}
Given  $f\in S(\ps)_\beta$ of big and nef degree $\beta=[D]\in A_{d-1}(\ps)$
and $\sigma\in\Sigma_D$,
let
$J^\sigma_0(f)$ be the ideal  in $S(\ps)$
generated by $x_\rho\partial f/\partial x_\rho$, $\rho\in\Sigma(1)$
 and all $x_{\rho'}$ such that $\rho'\subset\sigma$, and
let $J^\sigma_1(f)$
be the ideal quotient $J^\sigma_0(f):(\prod_{\rho\not\subset\sigma}x_\rho)$.
Then we get the quotient rings
$R_0^\sigma(f)=S(\ps)/J_0^\sigma(f)$ and
$R_1^\sigma(f)=S(\ps)/J_1^\sigma(f)$
graded by the Chow group $A_{d-1}(\ps)$.
\end{defn}

As a special case of \cite[Theorem~2.11]{m3}, we have:

\begin{thm}\label{t:main} Let $X$ be an anticanonical semiample
nondegenerate hypersurface defined by $f\in S_\beta$
in a complete simplicial toric variety $\ps$.
Then there is a natural isomorphism
$$\bigoplus_{p,q}H^{p,q}(X)\cong\bigoplus_{p,q} T(X)_{p,q}\oplus\biggl(\bigoplus_{\sigma\in\Sigma_X}
T^\sigma(X)_{s,s}
\otimes R^\sigma_1(f)_{(q-s)\beta+\beta_1^\sigma}\biggr),$$
where $s=(p+q-d+\dim\sigma+1)/2$ and $\beta_1^\sigma=
\deg(\prod_{\rho_k\subset\sigma}x_k)$.
\end{thm}

By the next statement, we can immediately see that all the building blocks
$R^\sigma_1(f)_{(q-s)\beta+\beta_1^\sigma}$ of
the cohomology of  partial resolutions in Theorem~\ref{t:main}
are independent of the  resolution and
intrinsic to an ample Calabi-Yau hypersurface:

\begin{pr}\label{p:iso} \cite{m3} Let $X$ be a big and nef
nondegenerate hypersurface defined by $f\in S_\beta$
in a complete  toric variety $\ps$ with  the associated contraction map
$\pi:\ps@>>>\psx$.
If $f_\sigma\in S(V(\sigma))_{\beta^\sigma}$
denotes the polynomial defining the hypersurface $\pi(X)\cap V(\sigma)$
in the toric variety
$V(\sigma)\subset\psx$ corresponding to $\sigma\in\Sigma_X$,
then, there is a natural isomorphism induced by the pull-back:
$$H^{d(\sigma)-*,*-1}H^{d(\sigma)-1}(\pi(X)\cap\ts)\cong
R_1(f_{\sigma})_{*\beta^{\sigma}-\beta_0^{\sigma}}{\cong}
R^\sigma_1(f)_{*\beta-\beta_0+\beta_1^\sigma},$$
where $d(\sigma)=d-\dim\sigma$, $\ts\subset V(\sigma)$ is
the maximal torus, and
$\beta_0$ and $\beta_0^{\sigma}$
denote the anticanonical degrees on $\ps$ and $V(\sigma)$, respectively.
\end{pr}

Given a mirror pair $(X,X^*)$ of smooth
Calabi-Yau hypersurfaces in Batyrev's construction,
we expect that, for a pair of cones $\sigma$ and $\sigma^*$ over the
dual faces of the reflexive polytopes $\Delta^*$ and $\Delta$,
$T^{\sigma}(X)_{s,s}$
with $s=(p+q-d+\dim\sigma+1)/2$,
in $H^{p,q}(X)$ interchanges, by the mirror map
(the isomorphism which maps the quantum cohomology of one Calabi-Yau
hypersurface to the B-model chiral ring of the other one), with
 $R^{\sigma^*}_1(g)_{(p+q-\dim\sigma^*)\beta^*/2+\beta_1^{\sigma^*}}$
in $H^{d-1-p,q}(X^*)$ (note that $\dim\sigma^*=d-\dim\sigma+1$),
where
$g\in S(\pp_{\Sigma^*})_{\beta^*}$
 determines $X^*$.
For the 0-dimensional cones $\sigma$ and $\sigma^*$, the interchange goes
between the {\it polynomial part} $R_1(g)_{*\beta^*}$  of
 one smooth Calabi-Yau hypersurface
and the toric part of the cohomology of the other one.
This correspondence
was already confirmed by the construction of the generalized
monomial-divisor
mirror map in \cite{m3}. On the other hand, one can deduce that
the dimensions of these spaces coincide for the pair of 3-dimensional
smooth Calabi-Yau hypersurfaces, by using Remark~5.3 in \cite{m1}.
The correspondence between the toric and polynomial parts was discussed in
\cite{ck}.

Now, let us turn our attention to a mirror pair of semiample
singular Calabi-Yau hypersurfaces $Y$ and $Y^*$. We know that
their string cohomology should have the same dimension as the
usual cohomology of possible crepant smooth resolutions $X$ and
$X^*$, respectively. Moreover, the A-model and B-model chiral
rings on the string cohomology should be isomorphic for $X$ and
$X^*$, respectively. We also know that the polynomial $g$
represents the complex structure of the hypersurface $Y^*$ and its
resolution $X^*$,
 and, by mirror symmetry, $g$ should correspond
to the complexified K\"ahler class
of the mirror Calabi-Yau hypersurface.
Therefore, based on the mirror correspondence of smooth
Calabi-Yau hypersurfaces, we make the following prediction for the small quantum ring
presentation
on the string cohomology space:
\begin{equation}\label{e:conj}
QH_{\rm st}^{p,q}(Y)\cong
\hspace{-0.05in}
\bigoplus_{(\sigma,\sigma^*)}
R_1(\omega_{\sigma^*})_{(p+q-\dim\sigma^*+2)\beta^{\sigma^*}/2-\beta_0^{\sigma^*}}\otimes
R_1(f_{\sigma})_{(q-p+d-\dim\sigma+1)\beta^\sigma/2-\beta_0^{\sigma}},
\end{equation}
where the sum is by all
pairs of the cones $\sigma$ and $\sigma^*$ (including 0-dimensional cones)
over the
dual faces of the reflexive polytopes, and where
$\omega_{\sigma^*}\in S(V(\sigma^*))_{\beta^{\sigma^*}}$ is a formal
restriction of $\omega\in S(\pp_{\Delta^*})_{\beta^*}$, which should be related to
the complexified K\"ahler class of the mirror (we will discuss this in Section~\ref{s:vs}).
This construction can be rewritten in simpler terms, which will help us
to give a conjectural description of the usual string cohomology space for all
semiample Calabi-Yau hypersurfaces.

First, recall Batyrev's presentation of the toric variety $\pp_\Delta$
for an {\it arbitrary} polytope $\Delta$ in $M$
(see \cite{b1}, \cite{c2}).
Consider the {\it Gorenstein} cone $K$ over
$\Delta\times\{1\}\subset M\oplus{\Bbb Z}$.
Let $S_\Delta$
be the subring of
${\Bbb C}[t_0,t_1^{\pm1},\dots,t_d^{\pm1}]$ spanned over $\Bbb C$ by all
monomials of the form $t_0^k t^m=t_0^kt_1^{m_1}\cdots t_d^{m_d}$ where $k\ge0$ and
$m\in k\Delta$. This ring is graded by the assignment
$\deg(t_0^k t^m)=k$.
Since the vector $(m,k)\in K$ if and only if $m\in k\Delta$,
the ring $S_\Delta$ is isomorphic to the semigroup algebra ${\Bbb C}[K]$.
The toric variety
 $\pp_\Delta$
can be represented as
$${\rm Proj}(S_\Delta)={\rm Proj}({\Bbb C}[K]).$$
The ring $S_\Delta$ has a nice connection to the homogeneous coordinate
ring
$S(\pp_\Delta)={\Bbb C}[x_\rho:\rho\in\Sigma_\Delta(1)]$
of the toric variety $\pp_\Delta$, corresponding to a fan $\Sigma_\Delta$.
If $\beta\in A_{d-1}(\pp_\Delta)$ is
the class of the  ample divisor
$\sum_{\rho\in\Sigma_\Delta(1)} b_\rho D_\rho$
 giving rise to
the polytope $\Delta$,
then
there is a natural isomorphism of graded rings
\begin{equation}\label{e:isom}
{\Bbb C}[K]\cong S_\Delta\cong\bigoplus_{k=0}^\infty S(\pp_\Delta)_{k\beta},
\end{equation}
sending $(m,k)\in\CC[K]_k$ to $t_0^k t^m$ and
$\prod_\rho x_\rho^{k b_\rho+\langle m,e_\rho\rangle}$,
where  $e_\rho$ is the minimal integral generator of the ray $\rho$.
Now, given $f\in S(\pp_\Delta)_{\beta}$, we get the ring
$R_1(f)$. The polynomial $f=\sum_{m\in\Delta}f(m)
x_\rho^{b_\rho+\langle m,e_\rho\rangle}$, where $f(m)$ are the coefficients,
 corresponds by the isomorphisms
(\ref{e:isom}) to
 $\sum_{m\in\Delta}f(m)t_0t^m\in (S_\Delta)_1$ and
$\sum_{m\in\Delta}f(m)[m,1]\in{\Bbb C}[K]_1$ (the brackets [{ }] are used to distinguish the lattice points from the vectors over $\CC$),
which we also denote by $f$.
By the proof of \cite[Theorem~11.5]{bc}, we have that
$$(S(\pp_\Delta)/J_0(f))_{k\beta}\cong
(S_\Delta/\langle t_i\partial f/\partial t_i:\, i=0,\dots,d\rangle)_k
\cong R_0(f,K)_k,$$
where $R_0(f,K)$ is the quotient of ${\Bbb C}[K]$ by the ideal generated
by  all ``logarithmic derivatives'' of $f$:
$$\sum_{m\in\Delta}((m,1)\cdot n) f(m)[m,1]$$ for $n\in N\oplus{\Bbb Z}$.
The  isomorphisms (\ref{e:isom}) induce the bijections
$$S(\pp_\Delta)_{k\beta-\beta_0}@>\prod_\rho x_\rho>>\langle
\prod_\rho x_\rho \rangle_{k\beta}\cong
(I_\Delta^{(1)})_k \cong {\Bbb C}[K^\circ]_k$$
($\beta_0=\deg(\prod_\rho x_\rho)$),
where $I_\Delta^{(1)}\subset S_\Delta$ is the ideal spanned by all monomials
$t_0^k t^m$ such that
$m$ is in the interior of $k\Delta$, and ${\Bbb C}[K^\circ]\subset{\Bbb C}[K]$
is the ideal
spanned by all lattice points in the relative interior of $K$.
Since the space $R_1(f)_{k\beta-\beta_0}$ is isomorphic to the image
of $\langle\prod_\rho x_\rho \rangle_{k\beta}$ in $(S(\pp_\Delta)/J_0(f))_{k\beta}$,
$$R_1(f)_{k\beta-\beta_0}\cong R_1(f,K)_k,$$
where $R_1(f,K)$ is the image of ${\Bbb C}[K^\circ]$
in the graded ring $R_0(f,K)$.

The above discussion applies well to all faces $\Gamma$ in $\Delta$.
In particular, if the toric variety $V(\sigma)\subset\pp_\Delta$ corresponds
to $\Gamma$, and $\beta^\sigma\in A_{d-\dim\sigma-1}(V(\sigma))$
is the restriction
of the ample class $\beta$, then
$$S(V(\sigma))_{*\beta^\sigma}\cong {\Bbb C}[C],$$
where $C$ is the Gorenstein cone over the polytope $\Gamma\times\{1\}$.
This induces an isomorphism
$$R_1(f_\sigma)_{*\beta^\sigma-\beta_0^\sigma}\cong R_1(f_C,C),$$
where $f_C=\sum_{m\in\Gamma} f(m)[m,1]$
in ${\Bbb C}[C]_1$ is the projection of $f$ to the cone $C$.

Now, we can restate our conjecture (\ref{e:conj}) in terms of Gorenstein
cones:
$$
\bigoplus_{p,q}QH_{\rm st}^{p,q}(Y)\cong\bigoplus_{\begin{Sb}
p,q\\
(C,C^*)\end{Sb}}
R_1(\omega_{C^*},C^*)_{(p+q-d+\dim C^*+1)/2}\otimes
R_1(f_{C},C)_{(q-p+\dim C)/2},
$$
where the sum is by all dual faces of the reflexive
Gorenstein cones $K$ and $K^*$.
This formula is already supported by Theorem~8.2 in \cite{bd},
which for ample Calabi-Yau hypersurfaces in weighted projective spaces
gives a corresponding decomposition of the stringy Hodge
numbers (see Remark~\ref{r:corrw} in the next section).
A generalization of \cite[Theorem~8.2]{bd} will be proved
in Section~\ref{section.bbo}, justifying the above conjecture
in the case of ample Calabi-Yau hypersurfaces in Fano toric varieties.

It is known that the string cohomology, which should be the limit of
the quantum cohomology ring, of smooth Calabi-Yau hypersurfaces
should be
 the same as the usual cohomology. We also know the property
that the quantum cohomology spaces should be isomorphic for the ample Calabi-Yau hypersurface
$Y$ and its crepant resolution $X$. Therefore, it makes sense to compare
the above description of $QH_{\rm st}^{p,q}(Y)$
 with the description of the cohomology of semiample Calabi-Yau hypersurfaces $X$
in Theorem~\ref{t:main}. We can see that
the right components in the tensor products coincide, by
Proposition~\ref{p:iso} and the definition of  $R_1(f_{C},C)$.
On the other hand, the left components in $QH_{\rm st}^{p,q}(Y)$
for the ample Calabi-Yau hypersurface
$Y$ do not depend on a resolution, while the left components $T^\sigma(X)$ in
 $H^{p,q}(X)$
for the resolution $X$ depend on the Stanley-Reisner ideal $SR(\Sigma)$.
This hints us to the following definitions:

\begin{defn}\label{d:rings}
Let $C$ be a Gorenstein cone in a lattice $L$,
subdivided by a fan $\Sigma$,
and let ${\Bbb C}[C]$ and ${\Bbb C}[C^\circ]$, where $C^\circ$ is the relative
interior of $C$,
 be the semigroup rings.
Define  ``deformed'' ring structures
$\CC[C]^\Sigma$ and $\CC[C]^\Sigma$ on ${\Bbb C}[C]$ and
${\Bbb C}[C^\circ]$, respectively, by the rule:
$[m_1][m_2]=[m_1+m_2]$ if  $m_1,m_2\subset\sigma\in\Sigma$, and
$[m_1][m_2]=0$, otherwise.

Given
 $g=\sum_{m\in C,\deg m=1} g(m)[m]$, where $g(m)$ are the coefficients, let
$$R_0(g,C)^\Sigma=\CC[C]^\Sigma/Z\cdot\CC[C]^\Sigma$$
be the graded ring over the graded module
$$R_0(g,C^\circ)^\Sigma=\CC[C^\circ]^\Sigma/Z\cdot\CC[C^\circ]^\Sigma,$$
where
$Z=\{\sum_{m\in C,\deg m=1} (m\cdot n)
g(m)[m]:\,n\in {\rm Hom}(L,{\Bbb Z})\}$.
Then define $R_1(g,C)^\Sigma$ as the image of
the natural homomorphism
$R_0(g,C^\circ)^\Sigma@>>>R_0(g,C)^\Sigma$.
\end{defn}

\begin{rem}
In the above definition, note that if $\Sigma$ is a trivial subdivision,
we recover the spaces $R_0(g,C)$ and $R_1(g,C)$ introduced earlier.
Also, we should mention that the Stanley-Reisner ring of the fan $\Sigma$ can be
naturally embedded into the ``deformed'' ring $\CC[C]^\Sigma$, and this map is an isomorphism
when the fan $\Sigma$ is smooth.
\end{rem}

Here is our conjecture about the string cohomology space of
semiample Calabi-Yau hypersurfaces in a complete toric variety.

\begin{conj}\label{semiampleconj}
Let $X\subset\ps$ be a semiample anticanonical nondegenerate hypersurface
defined by $f\in H^0(\ps,{\cal O}_\ps(X))\cong\CC[K]_1$,
and let $\omega$ be a generic element in $\CC[K^*]_1$, where $K^*$ is
the reflexive Gorenstein cone dual to the cone $K$ over the reflexive
polytope $\Delta$ associated to $X$.
Then there is a natural isomorphism:
$$
H^{p,q}_{\st}(X)\cong
\bigoplus_{C\subseteq K} R_1(\omega_{C^*},C^*)^\Sigma_{(p+q-d+\dim C^*+1)/2}
\otimes R_1(f_C,C)_{(q-p+\dim C)/2},
$$
where $C^*\subseteq K^*$ is a face dual to $C$, and
where $f_C$, $\omega_{C^*}$ denote the projections of $f$ and $\omega$
to the respective cones $C$ and $C^*$.
(Here, the superscript $\Sigma$ denotes
the subdivision of $K^*$ induced by the fan $\Sigma$.)
\end{conj}

Since the  dimension of the string cohomology for all crepant
partial resolutions should remain the same and should coincide
with the dimension of the quantum string cohomology space, we
expect that
\begin{equation}\label{e:expe}
{\dim}R_1(\omega_{C^*},C^*)_{\_}^\Sigma=
{\dim}R_1(\omega_{C^*},C^*)_{\_},
\end{equation}
which will be shown in  Section~\ref{section.brel} for a
projective subdivision $\Sigma$. Conjecture~\ref{semiampleconj}
will be confirmed by the corresponding decomposition of the
stringy Hodge numbers in Section~\ref{section.bbo}. Moreover, in
Section~\ref{s:vs},
 we will derive the Chen-Ruan orbifold cohomology as a special case
of Conjecture~\ref{semiampleconj} for ample Calabi-Yau hypersurfaces in complete simplicial
toric varieties.

\section{Hodge-Deligne numbers of affine hypersurfaces}\label{s:hd}

Here, we  compute the  dimensions of the spaces $R_1(g,C)_{\_}$
from the previous section.
It follows from Proposition~\ref{p:iso} that these dimensions are exactly
the Hodge-Deligne numbers of the minimal weight space on
the middle cohomology of a hypersurface in a torus.
An explicit formula in \cite{dk} and \cite{bd} for the $E$-polynomial
of a nondegenerate affine hypersurface whose Newton polyhedra is a simplex
leads us to the answer for the graded dimension of $R_1(g,C)$ when
$C$ is a simplicial Gorenstein cone.
However, it was very difficult to compute the Hodge-Deligne numbers
of  an arbitrary nondegenerate affine hypersurface. This was a major technical
problem in the proof of mirror symmetry of the stringy Hodge
numbers for  Calabi-Yau complete intersections in \cite{bb}.
Here, we will present a  simple formula for the Hodge-Deligne numbers
of a nondegenerate affine hypersurface.

Before we start computing ${\rm gr.dim.}R_1(g,C)$, let us note
that for a nondegenerate $g\in\CC[C]_1$ (i.e., the corresponding hypersurface
in ${\rm Proj}(\CC[C])$ is nondegenerate):
$${\rm gr.dim.}R_0(g,C)=S(C,t),$$
where the polynomial $S$ is the same as in Definition~\ref{d:bd}
of the stringy Hodge numbers. This was shown in
\cite[Theorem~4.8 and 2.11]{b1} (see also \cite{Bor.locstring}).

When the cone $C$ is simplicial, we already know the formula for the
graded dimension of $R_1(g,C)$:

\begin{pr}\label{p:simp} Let $C$ be a simplicial Gorenstein cone, and let
$g\in\CC[C]_1$ be nondegenerate.
Then
$${\rm gr.dim.}R_1(g,C)=\tilde S(C,t)$$
where
$\tilde S(C,t)=\sum_{C_1\subseteq C} S(C_1,t)
(-1)^{\dim C-\dim C_1}$.
\end{pr}

\begin{pf}
The polynomial $\tilde S(C,t)$ was introduced with a slightly
different notation
in \cite[Definition~8.1]{bd} for a lattice simplex.
One can  check that
$\tilde S(C,t)$ in this proposition is  equivalent to the one in
\cite[Corollary~6.6]{bd}.
{From} the previous section and \cite[Proposition~9.2]{b1}, we
know that
$$R_1(g,C)\cong\gr_F W_{\dim Z_g}H^{\dim Z_g}(Z_g),$$
where $Z_g$ is the nondegenerate affine hypersurface determined
by $g$ in the
maximal torus of  ${\rm Proj}(\CC[C])$.
By \cite[Proposition~8.3]{bd},
 $$E(Z_g;u,v)
=\frac{(uv-1)^{\dim C-1}+(-1)^{\dim C}}{uv}+(-1)^{\dim C}
\sum_{\begin{Sb}
C_1\subseteq C\\
\dim C_1>1\end{Sb}}\frac{u^{\dim C_1}}{uv}\tilde S(C_1,u^{-1}v).$$
Now, note that the coefficients $e^{p,q}(Z_g)$ at the monomials $u^p v^q$
with $p+q=\dim Z_g$ are related to the Hodge-Deligne
numbers by the calculations in \cite{dk}:
$$e^{p,q}(Z_g)=(-1)^{\dim C}h^{p,q}(H^{\dim Z_g}(Z_g))+(-1)^{p}\delta_{pq}
C_{\dim C-1}^p,$$
where $\delta_{pq}$ is the Kronecker symbol and  $C_{\dim C-1}^p$ is the
binomial coefficient.
Comparing this with the above formula for $E(Z_g;u,v)$,
we deduce the result.
\end{pf}

\begin{rem}\label{r:corrw}
By the above proposition, we can see that
\cite[Theorem~8.2]{bd} gives a decomposition of the
 stringy Hodge numbers of
ample Calabi-Yau hypersurfaces in weighted projective spaces
in  correspondence with Conjecture~\ref{semiampleconj}.
\end{rem}

Next, we generalize the polynomials $\tilde S(C,t)$
from Proposition~\ref{p:simp} to nonsimplicial Gorenstein cones in such a way
that
they would count the graded dimension of $R_1(g,C)$.

\begin{defn}\label{d:spol}

Let $C$ be a Gorenstein cone in a lattice $L$. Then set $$ \tilde
S(C,t) := \sum_{C_1\subseteq C} S(C_1,t) (-1)^{\dim C-\dim C_1}
G([C_1,C],t),$$ where   $G$ is a polynomial (from
Definition~\ref{Gpoly} in the Appendix) for the partially ordered
set $[C_1,C]$ of the faces of $C$ that contain $C_1$.
\end{defn}

\begin{rem}\label{tildepoincare}
It is not hard to show that the polynomial $\tilde S(C,t)$ satisfies
the duality
$$
\tilde S(C,t) = t^{\dim C} \tilde S(C,t^{-1})
$$
based on the duality properties of $S$ and the definition of $G$-polynomials.
However, the next result and Proposition~\ref{p:iso}
imply this fact.
\end{rem}

\begin{pr}\label{p:nonsimp} Let $C$ be a  Gorenstein cone, and let
$g\in\CC[C]_1$ be nondegenerate.
Then
$${\rm gr.dim.}R_1(g,C)=\tilde S(C,t).$$
\end{pr}

\begin{pf}
As in the proof of Proposition~\ref{p:simp},
we consider a nondegenerate affine hypersurface $Z_g$ determined
by $g$ in the
maximal torus of  ${\rm Proj}(\CC[C])$.
Then  \cite[Theorem 3.18]{bb} together with the definition of $S$ gives
$$E(Z_g;u,v)
= \frac{(uv-1)^{\dim C-1}}{uv} + \frac{(-1)^{\dim C}}{uv}
\sum_{C_2\subseteq C}
B([C_2,C]^*; u,v)S(C_2,vu^{-1})u^{\dim C_2},
$$
where the polynomials $B$ are from Definition~\ref{Q}.
 We use Lemma \ref{BfromG} and Definition \ref{d:spol} to rewrite this
as
\begin{multline*}
E(Z_g;u,v)
= \frac{(uv-1)^{\dim C-1}}{uv}
+\frac{(-1)^{\dim C}}{uv}\times
\\
\times
\sum_{C_2\subseteq C_1\subseteq C} u^{\dim C_2}
S(C_2,u^{-1}v)
 G([C_2,C_1],u^{-1}v)(-u)^{\dim C_1-\dim C_2}
G([C_1,C]^*,uv)
\\
=\frac{(uv-1)^{\dim C-1}}{uv}+\frac{(-1)^{\dim C}}{uv}
\sum_{C_1\subseteq C}u^{\dim C_1}
\tilde S(C_1,u^{-1}v)
G([C_1,C]^*,uv).
\end{multline*}

The definition of $G$-polynomials assures that
the degree of $u^{\dim C_1}G([C_1,C]^*,uv)$ is at most
$\dim C$ with the equality only when $C_1=C$.  Therefore,
the graded dimension of $R_1(g,C)$ can be read off the same way as
in the proof of Proposition~\ref{p:simp} from the
coefficients at total degree $\dim C-2$ in the above sum.
\end{pf}

\section{``Deformed''  rings and modules}
\label{section.brel}

While  this section may serve as an invitation to a new theory of
``deformed'' rings and modules, the goal here is to prove the
equality (\ref{e:expe}), by showing that the graded dimension
formula of Proposition~\ref{p:nonsimp} holds for the spaces
$R_1(g,C)^\Sigma$ from Definition~\ref{d:rings}. To prove the
formula we use the recent work of Bressler and Lunts (see
\cite{bl}, and also \cite{bbfk}). This requires us to first study
Cohen-Macaulay modules over the deformed semigroup rings
$\CC[C]^\Sigma$.

First, we want to generalize the nondegeneracy notion:

\begin{defn}\label{d:nond} Let $C$ be a Gorenstein cone in a lattice $L$,
subdivided by a fan $\Sigma$.
Given
 $g=\sum_{m\in C,\deg m=1} g(m)[m]$, get
$$g_j = \sum_{m\in C,\deg m=1}(m\cdot n_j) g(m)[m],\quad\text{ for }
j=1,\dots,\dim C$$
where $\{n_1,\dots,n_{\dim C}\}\subset{\rm Hom}(L,{\Bbb Z})$
descends to a basis of ${\rm Hom}(L,{\Bbb Z})/C^\perp$.
The element $g$ is called {\it $\Sigma$-regular (nondegenerate)} if
$\{g_1,\ldots,g_{\dim C}\}$ forms a regular sequence in the deformed semigroup ring
$\CC[C]^\Sigma$.
\end{defn}

\begin{rem}\label{r:nond} When $\Sigma$ is a trivial subdivision,
\cite[Theorem~4.8]{b1} shows that the above definition is consistent with
the  previous notion of nondegeneracy corresponding to the transversality
of a  hypersurface to torus orbits.
\end{rem}

\begin{thm}
\label{t:cm}
{\rm (i)} The  ring $\CC[C]^\Sigma$ and its module
$\CC[C^\circ]^\Sigma$ are Cohen-Macaulay.\\
{\rm (ii)} A generic element $g\in\CC[C]_1$ is $\Sigma$-regular. Moreover,
for a generic $g$ the sequence
$\{g_1,\ldots,g_{\dim C}\}$
from Definition~\ref{d:nond} is
$\CC[C^\circ]^\Sigma$-regular. \\
{\rm (iii)} If $g\in\CC[C]_1$ is $\Sigma$-regular, then the sequence
$\{g_1,\ldots,g_{\dim C}\}$
 is $\CC[C^\circ]^\Sigma$-regular.
\end{thm}

\begin{pf}

Part {\rm (ii)} follows from the proofs of Propositions 3.1 and 3.2
in \cite{Bor.locstring}.
The reader should notice that the proofs use degenerations defined by
projective simplicial subdivisions, and any fan
admits such a  subdivision.

Then, part  {\rm (ii)} implies {\rm (i)}, by the definition of Cohen-Macaulay,
while part {\rm (iii)} follows from {\rm (i)} and Proposition 21.9 in \cite{e}.
\end{pf}

As a corollary  of Theorem~\ref{t:cm},
we get the following simple description
of $\Sigma$-regular elements:

\begin{lem}\label{l:iff}
 An element $g\in\CC[C]_1$ is  $\Sigma$-regular, if and
only if its restriction to all maximum-dimensional cones $C'\in\Sigma(\dim C)$
is nondegenerate in $\CC[C']$.
\end{lem}

\begin{pf}
Since $\CC[C]^\Sigma$ is Cohen-Macaulay, the regularity of a sequence
is equivalent to the quotient by the sequence having a finite dimension, by
\cite[Theorem~17.4]{ma}.

One can check that $\CC[C]^\Sigma$ is filtered by the modules
$R_k$ defined as the span of $[m]$ such that the minimum cone that
contains $m$ has dimension at least $k$. The $k$-th graded
quotient of this filtration is the direct sum of $\CC[C_1^\circ]$
by all $k$-dimensional cones $C_1$ of $\Sigma$. If $g$ is
nondegenerate for every cone of maximum dimension, then its
projection to any cone $C_1$ is nondegenerate, and
Theorem~\ref{t:cm}
 shows that it is nondegenerate for each $\CC[C_1^\circ]$. Then
by decreasing induction on $k$ one shows that
$R_k/\{g_1,\ldots,g_{\dim C}\}R_k$ is finite-dimensional.

In the other direction, it is easy to see that for every $C'\in \Sigma$
the $\CC[C]^\Sigma$-module
$\CC[C']$ is a quotient of $\CC[C]^\Sigma$, which gives a surjection
$$\CC[C]^\Sigma/\{g_1,\ldots,g_{\dim C}\}\CC[C]^\Sigma
@>>>\CC[C']/\{g_1|_{C'},\dots,g_{\dim C}|_{C'}\}\CC[C']@>>>0.$$

\end{pf}

The above lemma implies that the property of $\Sigma$-regularity is preserved
by the restrictions:

\begin{lem} Let $C$ be a Gorenstein cone in a lattice $L$,
subdivided by a fan $\Sigma$.
 If $g\in\CC[C]_1$ is  $\Sigma$-regular, then
$g\in\CC[C_1]_1$ is $\Sigma$-regular for all faces $C_1\subseteq C$.
\end{lem}

\begin{pf}  Let $g\in\CC[C]_1$ be $\Sigma$-regular.
By Lemma~\ref{l:iff},
the restriction $g_{C'}$
is nondegenerate in $\CC[C']$ for all $C'\in\Sigma(\dim C)$.
Since the property of nondegeneracy associated with a hypersurface
is preserved by the restrictions,
$g_{C_1'}$
is nondegenerate in $\CC[C_1']$ for all  $C_1'\in\Sigma(\dim C_1)$ contained
in $C_1$. Applying Lemma~\ref{l:iff} again, we deduce the result.
\end{pf}

The next result generalizes \cite[Proposition~9.4]{b1} and
\cite[Proposition~3.6]{Bor.locstring}.

\begin{pr}\label{Zreg}
Let
$g\in\CC[C]_1$ be $\Sigma$-regular, then
$R_0(g,C)^\Sigma$ and $R_0(g,C^\circ)^\Sigma$
 have graded dimensions
$S(C,t)$ and $t^kS(C,t^{-1})$, respectively, and
there exists a nondegenerate pairing
$$
\langle\_,\_\rangle:
R_0(g,C)_k^\Sigma\times R_0(g,C^\circ)_{\dim C-k}^\Sigma\to R_0(g,C^\circ)_{\dim C}^\Sigma\cong\CC,
$$
induced by the multiplicative $R_0(g,C)^\Sigma$-module structure.
\end{pr}

\begin{pf}
It is easy to see that the above statement is equivalent to saying that
$\CC[C^\circ]^\Sigma$ is the canonical module for $\CC[C]^\Sigma$.
When $\Sigma$ consists of the faces of $C$ only, this is well-known
(cf. \cite{d}). To deal with the general case, we will heavily use
the results of \cite{e}, Chapter 21.

We denote $A=\CC[C]^\Sigma$. For every cone $C_1$ of $\Sigma$ the vector
spaces $\CC[C_1]$ and $\CC[C_1^\circ]$ are equipped with the natural $A$-module
structures. By Proposition 21.10 of \cite{e}, modified for the graded case,
we get
$$
\Ext^i_A(\CC[C_1],w_A) \iso \Bigl\{
\begin{array}{ll}
\CC[C_1^\circ],&i=\codim(C_1)\\
0,&i\neq\codim(C_1)
\end{array}
\Bigr.
$$
where $w_A$ is the canonical module of $A$.

Consider now the complex $\cal F$ of $A$-modules $$
0@>>>F^0@>>>F^1@>>>\cdots @>>>F^d@>>>0 $$ where $$ F^n =
\bigoplus_{C_1\in\Sigma,\codim(C_1)=n} \CC[C_1] $$ and the
differential is a sum of the restriction maps with signs according
to the orientations. The nontrivial cohomology of $\cal F$ is
located at $F^0$ and equals $\CC[C^\circ]^\Sigma$. Indeed, by
looking at each graded piece separately, we see that the
cohomology occurs only at $F^0$, and then the kernel of the map to
$F^1$ is easy to describe. We can now use the complex $\cal F$ and
the description of $\Ext^i_A(\CC[C_1],w_A)$ to try to calculate
$\Hom_A(\CC[C^\circ]^\Sigma,w_A)$. The resulting spectral sequence
degenerates immediately, and we conclude that
$\Hom_A(\CC[C^\circ]^\Sigma,w_A)$ has a filtration such that the
associated graded module is naturally isomorphic to $$
\bigoplus_{C_1\in \Sigma} \CC[C_1^\circ]. $$

By duality of maximal
Cohen-Macaulay modules (see \cite{e}), it suffices to show that
$\Hom_A(\CC[C^\circ]^\Sigma,w_A)\iso A$, but the above filtration
only establishes that it has the correct graded pieces, so extra
arguments are required.
Let $C^\prime$ be a cone of $\Sigma$ of maximum dimension. We observe that
$\cal F$ contains a subcomplex ${\cal F}^\prime$ such that
$$F^{\prime n}=
\bigoplus_{C_1\subseteq C^\prime} \CC[C_1].
$$
Similar to the case of $\cal F$, the cohomology of ${\cal F}^\prime$
occurs only at $F^{\prime 0}$ and equals $\CC[C^{\prime\circ}]$.
By snake lemma, the cohomology of
${\cal F}/{\cal F}^\prime$ also occurs at the zeroth spot and equals
$\CC[C^\circ]^\Sigma/\CC[C^{\prime\circ}]$. By looking at the spectral
sequences again, we see that
$$\Ext^{>0}(\CC[C^\circ]^\Sigma/\CC[C^{\prime\circ}],w_A)=0$$
and we have a grading preserving surjection
$$
\Hom_A(\CC[C^\circ]^\Sigma,w_A)
@>>>\Hom_A(\CC[C^{\prime\circ}],w_A)@>>>0.
$$
Since $\Hom_A(\CC[C^\prime],w_A)\iso \CC[C^{\prime\circ}]$,
duality of maximal Cohen-Macaulay modules over $A$ shows
that
$$\Hom_A(\CC[C^{\prime\circ}],w_A)\iso \CC[C^{\prime}]$$
so for every $m\in C^\prime$ the element $[m]$ of $A$ does not
annihilate the degree zero element of $\Hom_A(\CC[C^\circ]^\Sigma,w_A)$.
By looking at all $C^\prime$ together, this shows that
$$\Hom_A(\CC[C^\circ]^\Sigma,w_A)\iso A$$
which finishes the proof.
\end{pf}

\begin{pr}\label{p:nons} Let
$g\in\CC[C]_1$ be $\Sigma$-regular, then
the pairing $\langle\_\,,\_\rangle$
induces a symmetric nondegenerate pairing
$\{\_\,,\_\}$
on $R_1(g,C)^\Sigma$,
defined by
$$\{ x,y\} = \langle x,y'\rangle$$
where $y'$ is an element of $R_0(g,C^\circ)^\Sigma$ that maps to $y$.
\end{pr}

\begin{pf}
The nondegeneracy of the pairing $\{\_\,,\_\}$ follows from that of $\langle\_\,,\_\rangle$.
The pairing is symmetric, because it comes from the commutative product
on $\CC[C^\circ]^\Sigma$.
\end{pf}

\begin{thm}\label{dimtilde} Let $C$ be a  Gorenstein cone
subdivided by a projective fan $\Sigma$. If
$g\in\CC[C]_1$ is $\Sigma$-regular, then
the graded dimension of $R_1(g,C)^\Sigma$ is $\tilde S(C,t)$.
\end{thm}

\begin{pf} We will use the description of Bressler and Lunts
\cite{bl} of locally free flabby sheaves on the finite ringed
topological space associated to the cone $C$. We recall here
the basic definitions. Consider the set $P$ of all faces of the cone $C$.
It is equipped with the topology in which open sets are subfans,
i.e. the sets of faces closed under the operation of taking a face.
Bressler and Lunts define a sheaf $\cal A$ of graded commutative rings on $P$
whose sections over each open set is the ring of continuous
piecewise polynomial functions on the union of all strata of this set.
The grading of linear functions will be set to $1$, contrary
to the convention of \cite{bl}.

They further restrict their attention to the sheaves $\cal F$ of
$\cal A$-modules on $P$ that satisfy the following conditions.

\noindent$\bullet$
For every face $C_1$ of $C$, sections of $\cal F$ over the open set
that corresponds to the union of all faces of $C_1$ is
a free module over the ring of polynomial functions on $C_1$.

\noindent$\bullet$
${\cal F}$ is flabby, i.e. all restriction maps are surjective.

We will use the following crucial result.
\begin{thm}\cite{bl}
Every sheaf $\cal F$ that satisfies the above two properties is
isomorphic to a direct sum of indecomposable graded sheaves ${\cal
L}_{C_1}t^i$, where $C_1$ is a face of $C$ and $t^i$ indicates a
shift in grading. For each indecomposable sheaf $L_{C_1}$ the
space of global sections $\Gamma(P,{\cal L}_{C_1})$ is a module
over the polynomial functions on $C$ of the graded rank
$G([C_1,C]^*,t)$ where $[C_1,C]$ denotes the Eulerian subposet of
$P$ that consists of all faces of $C$ that contain $C_1$.
\end{thm}

Now let us define a sheaf ${\cal B}(g)$ on $P$ whose sections over the
open subset $I\in P$ are $\CC[\cap_{i\in I}C_i]^\Sigma$. It is clearly
a flabby sheaf, which can be given a grading by $\deg(\_)$.
Moreover, ${\cal B}(g)$ can be given a structure of a
sheaf of $\cal A$ modules as follows. Every linear function $\varphi$ on
a face $C_1$ defines a logarithmic derivative
$$
\partial_\varphi g:=\sum_{m\in C_1,\deg m =1} \phi(m)g(m)[m]
$$
of $g$, which is an element
of the degree $1$ in $\CC[C_1]^\Sigma$. Then the action of $\varphi$
is given by the multiplication by $\partial_\varphi g$, and this
action is extended to all polynomial functions on the cone $C_1$.
Similar construction clearly applies to continuous piecewise polynomial
functions for any open set of $P$.

Proposition~\ref{Zreg} assures that ${\cal B}(g)$ satisfies
the second condition of Bressler and Lunts, and can therefore be decomposed
into a direct sum of $L_{C_1}t^i$ for various $C_1$ and $i$.
The definition of $R_1(g,C)^\Sigma$ implies that its graded dimension is equal
to the graded rank of the stalk of ${\cal B}(g)$ at the point $C\in P$.
Since the graded rank of $\cal B$ is $S(C,t)$, we conclude that
$$
S(C,t) = \sum_{C_1\subseteq C} {\rm gr.dim.}R_1(g_{C_1},C_1)^\Sigma G([C_1,C]^*,t).
$$
To finish the proof of Theorem \ref{dimtilde}, it
remains to apply Lemma \ref{Ginverse}.
\end{pf}

\section{Decomposition of stringy Hodge numbers for hypersurfaces}
\label{section.bbo}

In this section, we prove a generalization of \cite[Theorem~8.2]{bd}
for all  Calabi-Yau hypersurfaces, which gives a
decomposition of the stringy Hodge numbers of the hypersurfaces.
First, we  recall a formula for the stringy Hodge numbers
of Calabi-Yau hypersurfaces obtained in \cite{bb}. Then using
a bit of combinatorics, we rewrite this formula precisely to
the form of \cite[Theorem~8.2]{bd} with $\tilde S$ defined in the previous section.

 The stringy
Hodge numbers of a Calabi-Yau complete intersection
 have been calculated in \cite{bb} in
terms of the numbers of integer points inside multiples of various
faces of the reflexive polytopes
$\Delta$ and $\Delta^*$ as well as some polynomial invariants of
partially ordered sets.
A special case of the main result in \cite{bb} is the following description
of the stringy $E$-polynomials of  Calabi-Yau hypersurfaces.

\begin{thm}\label{st.formula} \cite{bb} Let $K\subset M\oplus{\Bbb Z}$ be
the Gorenstein cone
 over  a reflexive  polytope $\Delta\subset M$.
For every $(m,n)\in (K,K^*)$ with $m\cdot n =0$ denote by
$x(m)$ the minimum face of $K$ that contains $m$ and by $x^*(n)$
the dual of the minimum face of $K^*$ that contains $n$.
Also, set $A_{(m,n)}(u,v)$ be
$$
\frac{(-1)^{\dim(x^*(n))}}{uv}
(v-u)^{\dim(x(m))}(uv -1)^{d+1-\dim(x^*(n))}B([x(m),x^*(n)]^*;u,v)
$$
where the function  $B$ is defined
in Definition~\ref{Q} in the Appendix.
Then
$$
E_{\rm st}(Y; u,v)=
\sum_{(m,n) \in (K,K^*),m\cdot n =0}
\left(\frac{u}{v}\right)^{{\rm deg}\,m} A_{(m,n)}(u,v)
\left(\frac{1}{uv}\right)^{{\rm deg}\,n}
$$
for an  ample nondegenerate Calabi-Yau hypersurface $Y$ in
$\PP_\Delta={\rm Proj}(\CC[K])$.
\end{thm}

The mirror duality $E_{\st}(Y;u,v) = (-u)^{d-1}E_\st(Y^*;u^{-1},v)$
was proved in \cite{bb} as the immediate corollary of the above formula
and the duality property $B(P; u,v) = (-u)^{\rk P} B(P^*;u^{-1},v)$.
It was not noticed there that Lemma \ref{Ginverse} allows one to
rewrite the $B$-polynomials in terms of $G$-polynomials, which we will
now use to give a formula for the $E_\st(Y;u,v)$,  explicitly
obeying the mirror duality.
The next result is a generalization of
Theorem~8.2 in \cite{bd} with $\tilde S$ from Definition~\ref{d:spol}.

\begin{thm}\label{estfromtilde}
Let $Y$ be an ample nondegenerate Calabi-Yau hypersurface  in
$\PP_\Delta={\rm Proj}(\CC[K])$. Then
$$
E_\st(Y;u,v) = \sum_{C\subseteq K}
(uv)^{-1}(-u)^{\dim C}\tilde S(C,u^{-1}v)
\tilde S(C^*,uv).
$$
\end{thm}

{\em Proof.} First,
 observe that the formula for $E_{\rm st}(Y; u,v)$ from
 Theorem \ref{st.formula}
 can be written as
$$ \sum_{m,n,C_1,C_2} \frac{(-1)^{\dim C_2^*}}{uv}
(v-u)^{\dim(C_1)}B([C_2,C_1^*];u,v)(uv -1)^{\dim C_2}
\left(\frac{u}{v}\right)^{{\rm deg} m}
\left(\frac{1}{uv}\right)^{{\rm deg}n} $$ where the sum is taken
over all pairs of cones $C_1\subseteq K, C_2\subseteq K^*$ that
satisfy $C_1\cdot C_2 = 0$ and all $m$ and $n$ in the relative
interiors of $C_1$ and $C_2$, respectively. We use the standard
duality result (see Definition~\ref{d:bd})
 $$\sum_{n\in {\rm int}(C)}
t^{-\deg(n)}= (t-1)^{-\dim C}S(C,t)$$ to rewrite the above formula
as $$\frac 1{uv}\sum_{C_1\cdot C_2=0} (-1)^{\dim(C_2^*)}u^{\dim
C_1} B([C_2,C_1^*];u,v)S(C_1,u^{-1}v)S(C_2,uv). $$ Then apply
Lemma~\ref{BfromG} to get $$ E_\st(Y;u,v) = \frac 1{uv} \sum_{C\in
K} \sum_{C_1\subseteq C,C_2\subseteq C^*}
(-1)^{\dim(C_2^*)}u^{\dim C_1} ~\times $$ $$ \times~
G([C_1,C],u^{-1}v) (-u)^{\dim C_1^* -\dim C^*} G([C_2,C^*],uv)
S(C_1,u^{-1}v)S(C_2,uv). $$ It remains to use
Definition~\ref{d:spol}. \qed

\section{String cohomology construction via intersection cohomology}
\label{section.general}

Here, we construct the string cohomology space for
$\QQ$-Gorenstein  toroidal varieties, satisfying the assumption of
Proposition~\ref{BDvsB}.
 The motivation for this construction comes from the conjectural
description of the string cohomology space for
ample Calabi-Yau hypersurfaces and a look at the formula in
\cite[Theorem~6.10]{bd} for the stringy
$E$-polynomial of a Gorenstein variety with  abelian quotient singularities.
This immediately leads to a decomposition of the string cohomology space
as a direct sum of  tensor products of the usual cohomology of a closure
of a strata with the spaces $R_1(g,C)$ from Proposition~\ref{p:simp}.
Then the  property that the
 intersection cohomology of an orbifold is naturally isomorphic to
the usual cohomology leads us to the construction of the string
cohomology space for $\QQ$-Gorenstein toroidal varieties. We show
that this space has the dimension prescribed by
Definition~\ref{d:bd} for Gorenstein complete toric varieties and
the  nondegenerate complete intersections in them.

\begin{cdefn}\label{d:orbdef} Let $X=\bigcup_{i \in I} X_i$ be
a Gorenstein complete variety with  quotient abelian
singularities, satisfying the assumption of
Proposition~\ref{BDvsB}.
 The stringy Hodge spaces of $X$ are naturally isomorphic to
$$H_{\rm st}^{p,q}(X)\cong\bigoplus_{\begin{Sb}i \in I\\
k\ge0 \end{Sb}}H^{p-k,q-k}(\overline{X}_i)\otimes
R_1(\omega_{\sigma_i},\sigma_i)_k,$$
where $\sigma_i$ is the Gorenstein simplicial cone of the singularity along
the strata $X_i$, and $\omega_{\sigma_i}\in\CC[\sigma_i]_1$ are  nondegenerate such that,
for $\sigma_j\subset \sigma_i$,
$\omega_{\sigma_i}$ maps to $\omega_{\sigma_j}$ by the natural projection
$\CC[\sigma_i]@>>>\CC[\sigma_j]$.
\end{cdefn}

\begin{rem}
Since $\overline{X}_i$ is a compact orbifold, the coefficient
$e^{p,q}(\overline{X}_i)$ at the monomial  $u^p v^q$ in the polynomial
$E(\overline{X}_i;u,v)$ is equal to $(-1)^{p+q}h^{p,q}(\overline{X}_i)$,
by Remark~\ref{r:epol}. Therefore, Proposition~\ref{p:simp} shows that
the above decomposition of  $H_{\rm st}^{p,q}(X)$ is in correspondence
with \cite[Theorem~6.10]{bd}, and the dimensions $h_{\rm st}^{p,q}(X)$
coincide with
those from Definition~\ref{d:bcangor}.
\end{rem}

Since we expect that the usual cohomology must be replaced in
Definition~\ref{d:orbdef} by the intersection
cohomology for Gorenstein toroidal varieties, the next result is
a natural generalization of Theorem~6.10 in \cite{bd}.

\begin{thm} \label{8.3}
Let $X=\bigcup_{i \in I} X_i$ be a Gorenstein complete  toric
variety or a  nondegenerate complete intersection of Cartier
hypersurfaces in the toric variety, where the stratification is
induced by the torus orbits.
 Then
 $$E_{\rm st}(X;u,v)= \sum_{i \in I} E_{\rm int}(\overline X_i;u,v)
  \cdot \tilde S(\sigma_i,uv),$$
  where
$\sigma_i$ is the Gorenstein cone of the singularity along the
strata $X_i$.
\end{thm}

\begin{pf}
Similarly to Corollary~3.17 in  \cite{bb}, we have $$E_{\rm
int}(\overline X_i;u,v)=\sum_{X_j\subseteq \overline
X_i}E(X_i;u,v)\cdot G([\sigma_i\subseteq \sigma_j]^*,uv).$$
 Hence, we get
 $$ \sum_{i \in I} E_{\rm
int}(\overline X_i;u,v) \cdot \tilde S(\sigma_i,uv) = \sum_{i\in
I}\sum_{X_j\subseteq \overline X_i} E(X_j;u,v)
G([\sigma_i\subseteq \sigma_j]^*,uv)\tilde S(\sigma_i,uv) $$ $$
=\sum_{j\in I} E(X_j;u,v) \Bigl(\sum_{\sigma_i\subseteq \sigma_j}
G([\sigma_i\subseteq \sigma_j]^*,uv)\tilde S(\sigma_i,uv) \Bigr) =
\sum_{j\in I} E(X_j;u,v) S(\sigma_j,uv), $$ where at the last step
we have used the formula for $\tilde S$ and Lemma \ref{Ginverse}.
\end{pf}


Based on the above theorem, we propose the following conjectural
description of the stringy Hodge spaces for $\QQ$-Gorenstein
toroidal varieties.

\begin{cdefn}\label{d:tordef} Let $X=\bigcup_{i \in I} X_i$ be
 a $\QQ$-Gorenstein $d$-dimensional complete toroidal variety, satisfying the assumption of
Proposition~\ref{BDvsB}. The stringy Hodge spaces of $X$ are defined by:
$$H_{\rm st}^{p,q}(X):=\bigoplus_{\begin{Sb}i \in I\\
k\ge0 \end{Sb}}H_{\rm int}^{p-k,q-k}(\overline{X}_i)\otimes
R_1(\omega_{\sigma_i},\sigma_i)_k,$$
where $\sigma_i$ is the Gorenstein  cone of the singularity along
the strata $X_i$, and $\omega_{\sigma_i}\in\CC[\sigma_i]_1$ are  nondegenerate such that,
for $\sigma_j\subset \sigma_i$,
$\omega_{\sigma_i}$ maps to $\omega_{\sigma_j}$ by the natural projection
$\CC[\sigma_i]@>>>\CC[\sigma_j]$. Here, $p,q$ are rational numbers from $[0,d]$, and
we assume that $H_{\rm int}^{p-k,q-k}(\overline{X}_i)=0$ if $p-k$ or $q-k$ is not a
 non-negative integer.
\end{cdefn}

\begin{rem}
Toric varieties and nondegenerate complete intersections of Cartier hypersurfaces have
the stratification induced by the torus orbits which   satisfies
the assumptions in the above definition.
\end{rem}

\section{String cohomology vs. Chen-Ruan orbifold cohomology}\label{s:vs}

Our next goal is to compare the two descriptions of string cohomology for
 Calabi-Yau  hypersurfaces to the Chen-Ruan orbifold cohomology.
Using the work of \cite{p},
we will show that in the case of ample orbifold Calabi-Yau  hypersurfaces
the three descriptions coincide.
We refer the reader to \cite{cr} for the orbifold cohomology
theory and only use \cite{p} in order to describe the orbifold
cohomology for complete simplicial toric varieties and Calabi-Yau
hypersurfaces in Fano simplicial toric varieties.

From Theorem~1 in \cite[Section~4]{p} and the definition of the orbifold Dolbeault  cohomology
space we deduce:

\begin{pr}
Let $\ps$ be a $d$-dimensional complete simplicial toric variety.
Then the orbifold Dolbeault  cohomology space of $\ps$ is
$$H^{p,q}_{orb}(\ps;\CC)\cong\bigoplus_{\begin{Sb}\sigma\in\Sigma\\l\in\QQ\end{Sb}}
H^{p-l,q-l}(V(\sigma))\otimes \bigoplus_{t\in T(\sigma)_l}\CC t,$$
where $T(\sigma)_l=\{\sum_{\rho\subset\sigma}a_\rho [e_\rho]\in
N:a_\rho\in(0,1), \sum_{\rho\subset\sigma}a_\rho=l\}$ (when
$\sigma=0$, set $l=0$ and $T(\sigma)_l=\CC$), and $V(\sigma)$ is
the closure of the torus orbit corresponding to $\sigma\in\Sigma$.
 Here, $p$ and $q$ are rational numbers in $[0,d]$,
and $H^{p-l,q-l}(V(\sigma))=0$ if $p-l$ or $q-l$ is not integral. (The elements of
$\oplus_{0\ne\sigma\in\Sigma,l}T(\sigma)_l$ correspond to the twisted sectors.)
\end{pr}

In order to compare this result to the description in
Definition~\ref{d:tordef}, we need to specify the
$\omega_{\sigma_i}$ for the toric variety $\ps$. The
stratification of $\ps$ is given by the torus orbits:
$\ps=\cup_{\sigma\in\Sigma}\TT_\sigma$. The singularity of the
variety $\ps$ along the strata $\TT_\sigma$ is given by the cone
$\sigma$, so we need to specify a nondegenerate
$\omega_\sigma\in\CC[\sigma]_1$ for each $\sigma\in\Sigma$. If
$\omega_\sigma=\sum_{\rho\subset\sigma}\omega_\rho [e_\rho]$ with
$\omega_\rho\ne0$, then one can deduce that $\omega_\sigma$ is
nondegenerate using Remark~\ref{r:nond} and the fact that the
nondegeneracy of a  hypersurface in a complete simplicial toric
variety (in this case, it  corresponds to a simplex) is equivalent
to the nonvanishing of the logarithmic derivatives simultaneously.
So, picking any nonzero coefficients $\omega_\rho$ for each
$\rho\in\Sigma(1)$ gives a nondegenerate
$\omega_\sigma\in\CC[\sigma]_1$ satisfying the condition of
Definition~\ref{d:tordef}. For such $\omega_\sigma$, note that the
set $Z=\{\sum_{\rho\subset\sigma} (e_\rho\cdot m) \omega_\rho
e_\rho:\,m\in {\rm Hom}(N,{\Bbb Z})\}$ is a linear span of
$e_\rho$ for $\rho\subset\sigma$. Hence,
$$R_0(\omega_\sigma,\sigma)_l=(\CC[\sigma]/Z\cdot\CC[\sigma])_l\cong
\bigoplus_{t\in \tilde{T}(\sigma)_l}\CC t,$$ where $\tilde
T(\sigma)_l=\{\sum_{\rho\subset\sigma}a_\rho e_\rho\in
N:a_\rho\in[0,1), \sum_{\rho\subset\sigma}a_\rho=l\}$, and
$$R_1(\omega_\sigma,\sigma)_l\cong
\bigoplus_{t\in{T}(\sigma)_l}\CC t.$$ This shows that the orbifold
Dolbeault  cohomology for complete simplicial toric varieties can
be obtained as a special case of the description of string
cohomology in Definition~\ref{d:tordef}.

We will now explain how  the parameter $\omega$ should be related
to the complexified K\"ahler class.  We do not have the definition
of the "orbifold" K\"ahler cone even for simplicial toric
varieties. However, we know the K\"ahler classes in
$H^2(\ps,\RR)$.

\begin{pr}
Let $\ps$ be a projective simplicial toric variety,
then $H^2(\ps,\RR)\cong PL(\Sigma)/M_\RR$,
where   $PL(\Sigma)$ is the set of $\Sigma$-piecewise linear functions $\varphi:N_\RR@>>>\RR$,
which are linear on each $\sigma\in\Sigma$. The K\"ahler cone $K(\Sigma)\subset H^2(\ps,\RR)$
of $\ps$ consists of the classes of the upper strictly convex $\Sigma$-piecewise linear functions.
\end{pr}

One may call $K(\Sigma)$ the ``untwisted'' part of the orbifold
K\"ahler cone. So, we can introduce the {\it untwisted
complexified K\"ahler space} of the complete simplicial toric
variety: $$K^{\rm untwist}_\CC(\ps)= \{\omega\in
H^2(\ps,\CC):Im(\omega)\in K(\Sigma)\}/{\rm im}H^2(\ps,\ZZ).$$ Its
elements may be called the {\it untwisted complexified K\"ahler
classes}. We can find a generic enough $\omega\in K^{\rm
untwist}_\CC(\ps)$ represented by a complex valued
$\Sigma$-piecewise linear function $\varphi_\omega:N_\CC@>>>\CC$
such that $\varphi_\omega(e_\rho)\ne0$ for $\rho\in\Sigma(1)$.
Setting $\omega_\rho=\exp(\varphi_\omega(e_\rho))$ produces our
previous parameters $\omega_\sigma$ for $\sigma\in\Sigma$. This is
how we believe $\omega_\sigma$ should relate to the complexified
K\"ahler classes, up to perhaps some instanton corrections.

We next turn our attention to the case of an ample Calabi-Yau hypersurface $Y$ in a complete simplicial
toric variety $\ps$. Section 4.2 in \cite{p} works with a generic nondegenerate anticanonical
hypersurface. However, one can avoid the use of Bertini's theorem and state the result without
``generic''. It is shown that the  nondegenerate anticanonical
hypersurface $X$ is a suborbifold of $\ps$, the twisted sectors of $Y$ are obtained by intersecting
with the closures of the torus orbits and the degree shifting numbers are the same as for the
toric variety $\ps$. Therefore, we conclude:

\begin{pr}\label{p:des}
Let $Y\subset\ps$ be an ample Calabi-Yau hypersurface in a complete simplicial
toric variety.
Then
$$H^{p,q}_{orb}(Y;\CC)\cong\bigoplus_{\begin{Sb}\sigma\in\Sigma\\l\in\ZZ\end{Sb}}
H^{p-l,q-l}(Y\cap V(\sigma))\otimes
\bigoplus_{t\in T(\sigma)_l}\CC t,$$
where $T(\sigma)_l=\{\sum_{\rho\subset\sigma}a_\rho [e_\rho]\in N:a_\rho\in(0,1),
\sum_{\rho\subset\sigma}a_\rho=l\}$ (when $\sigma=0$, set $l=0$ and $T(\sigma)_l=\CC$).
\end{pr}

As in the case of the toric variety, we pick
$\omega_\sigma=\sum_{\rho\subset\sigma}\omega_\rho e_\rho$
with  $\omega_\rho\ne0$. Then, by the above proposition,
$$H_{\rm st}^{p,q}(Y)\cong H^{p,q}_{orb}(Y;\CC).$$

We now want to show that the description in Proposition~\ref{p:des}
 is equivalent to the one in Conjecture~\ref{semiampleconj}.
First, note that the proper faces $C^*$ of the Gorenstein cone $K^*$ in Conjecture~\ref{semiampleconj}
one to one correspond to the cones $\sigma\in\Sigma$. Moreover, the rings $\CC[C^*]\cong\CC[\sigma]$
are isomorphic
in this correspondence. If we take $\omega\in\CC[K^*]^\Sigma_1$
to be
$[0,1]+\sum_{\rho\in\Sigma(1)}\omega_\rho[e_\rho,1]$,
then $\omega$ is $\Sigma$-regular and
$$R_1(C^*,\omega_{C^*})_l\cong R_1(\omega_\sigma,\sigma)_l\cong
\oplus_{t\in{T}(\sigma)_l}\CC t.$$
On the other hand, the Hodge component $H^{p-l,q-l}(Y\cap V(\sigma))$ decomposes into the
direct sum
$$H^{p-l,q-l}_{\rm toric}(Y\cap V(\sigma))\oplus H^{p-l,q-l}_{\rm res}(Y\cap V(\sigma))$$
of the toric and residue parts.
Since $Y\cap V(\sigma)$ is an ample hypersurface, from \cite[Theorem~11.8]{bc} and
 Section~\ref{section.anvar} it follows  that
$$H^{p-l,q-l}_{\rm res}(Y\cap V(\sigma))\cong R_1(f_C,C)_{q-l+1},$$
where $C\subset K$ is the face dual to $C^*$ which corresponds to $\sigma$,
$p+q-2l=\dim Y\cap V(\sigma)=d-\dim\sigma-1=d-\dim C^*-1$.
If $p+q-2l\ne d-\dim C^*-1$, then $H^{p-l,q-l}_{\rm res}(Y\cap V(\sigma))=0$.
Hence, we get
$$\bigoplus_{\begin{Sb}\sigma\in\Sigma\\l\in\ZZ\end{Sb}}
H_{\rm res}^{p-l,q-l}(Y\cap V(\sigma))\otimes
\bigoplus_{t\in T(\sigma)_l}\CC t
\cong \bigoplus_{0\ne C\subseteq K} R_1(\omega_{C^*},C^*)^\Sigma_{a}
\otimes R_1(f_C,C)_{b},$$
where $a=(p+q-d+\dim C^*+1)/2$ and $b=(q-p+\dim C)/2$. We are left to show that
\begin{equation}\label{e:lef}
\bigoplus_{\begin{Sb}\sigma\in\Sigma\\l\in\ZZ\end{Sb}}
H_{\rm toric}^{p-l,p-l}(Y\cap V(\sigma))\otimes
\bigoplus_{t\in T(\sigma)_l}\CC t
\cong R_1(\omega,K^*)^\Sigma_{p+1}.
\end{equation}
Notice that the dimensions of the spaces on both sides coincide, so it suffices
to construct a surjective map between them.
This  will follow from the following proposition.

\begin{pr}\label{p:sttor} Let $\ps={\rm Proj}(\CC[K])$ be the Gorenstein Fano simplicial toric variety,
where $K$ as above.
Then there is a natural isomorphism:
$$H^{p,p}_{\rm st}(\ps)\cong\bigoplus_{\begin{Sb}\sigma\in\Sigma\\l\in\ZZ\end{Sb}}
H^{p-l,p-l}(V(\sigma))\otimes \bigoplus_{t\in T(\sigma)_l}\CC t\cong
R_0(\omega,K^*)_p^\Sigma,$$
where $\omega=[0,1]+\sum_{\rho\in\Sigma(1)}\omega_\rho[e_\rho,1]$ with $\omega_\rho\ne0$.
\end{pr}

\begin{pf} First, observe that the dimensions of the spaces in the isomorphisms coincide
by our definition of string cohomology, Proposition~\ref{Zreg} and \cite[Theorem~7.2]{bd}.
So, it suffices to construct a surjective map between them.

We know the cohomology ring of the toric variety:
$$H^{*}(V(\sigma))\cong
\CC[D_\rho:\rho\in\Sigma(1),\rho+\sigma\in\Sigma(\dim\sigma+1)]/
(P(V(\sigma))+SR(V(\sigma))),$$
where $$SR(V(\sigma))=\bigl\langle
D_{\rho_1}\cdots D_{\rho_k}:\{e_{\rho_1},\dots,e_{\rho_k}\}
\not\subset\tau \text{ for all
}\sigma\subset\tau\in\Sigma(\dim\sigma+1)\bigr\rangle$$ is the
Stanley-Reisner ideal, and $$P(V(\sigma))=\biggl\langle
\sum_{\rho\in\Sigma(1),\rho+\sigma\in\Sigma(\dim\sigma+1)} \langle
m,e_\rho\rangle D_\rho: m\in M\cap\sigma^\perp\biggr\rangle.$$

Define the maps from $H^{p-l,p-l}(V(\sigma))\otimes \bigoplus_{t\in T(\sigma)_l}\CC t$
to $R_0(\omega,K^*)^\Sigma$ by sending
 $D_{\rho_1}\cdots D_{\rho_{p-l}}\otimes t$ to
$\omega_{\rho_1}[e_{\rho_1}]\cdots \omega_{\rho_{p-l}} [e_{\rho_{p-l}}]\cdot t\in\CC[N]^\Sigma$.
One can easily see that these maps are well defined.
To finish the proof we need to show that the images cover $R_0(\omega,K^*)^\Sigma$.
Every lattice point $[n]$ in the boundary of $K^*$ lies in the relative interior of a face $C\subset K^*$,
and
can be written as a linear combination of the minimal integral generators of $C$:
$$[n]=\sum_{[e_\rho,1]\in C}(a_\rho+b_\rho)[e_\rho,1],$$
where  $a_\rho\in(0,1)$ and $b_\rho$ are nonnegative integers.
Let $C'\subseteq C$ be the cone spanned by those  $[e_\rho,1]$ for which
$a_\rho\ne0$.
The lattice point $\sum_{[e_\rho,1]\in C'}(a_\rho)[e_\rho,1]$ projects to one of the elements $t$
from $T(\sigma)_l$ for  some $l$ and $\sigma$ corresponding to $C'$.
Using the relations
$\sum_{\rho\in\Sigma(1)}\omega_\rho\langle m,e_\rho\rangle [e_\rho,1]$ in
the ring $R_0(\omega,K^*)^\Sigma$, we get that
$$[n]=\sum_{[e_\rho,1]\in C'}(a_\rho)[e_\rho,1]+
\sum_{\rho+\sigma\in\Sigma(\dim\sigma+1)}b'_\rho[e_\rho,1],$$
which comes from $H^{p-l,p-l}(V(\sigma))\otimes \CC t$
for an appropriate $p$.
The surjectivity now follows from the fact that
the boundary points of $K^*$ generate the ring
$C[K^*]^\Sigma/\langle \omega\rangle$.
\end{pf}

The isomorphism (\ref{e:lef}) follows from the above proposition
and the presentation:
$$H_{\rm toric}^{*}(Y\cap V(\sigma))\cong H^*(V(\sigma))/{\rm Ann}([Y\cap V(\sigma)])$$
(see (\ref{e:ann})).
Indeed, the map constructed in the proof of Proposition~\ref{p:sttor} produces a well defined
map between the right hand side
in (\ref{e:lef}) and $R_0(\omega,K^*)^\Sigma/Ann([0,1])$
 because the annihilator of $[Y\cap V(\sigma)]$ maps to the annihilator
of $[0,1]$. On the other hand,
$$(R_0(\omega,K^*)^\Sigma/Ann([0,1]))_p\cong R_1(\omega,K^*)^\Sigma_{p+1},$$
which is induced by the multiplication by $[0,1]$ in $R_0(\omega,K^*)^\Sigma$.

\begin{conj} We expect that the product structure on
$H^{*}_{\rm st}(\ps)$ is given by the ring structure $R_0(\omega,K^*)^\Sigma$.
Also, the ring structure on $R_1(\omega,K^*)^\Sigma_{*+1}$
induced from $R_0(\omega,K^*)^\Sigma/Ann([0,1])$ should give a subring of
$H^{*}_{\rm st}(Y)$ for a generic $\omega$ in Conjecture~\ref{semiampleconj}.
Moreover,
$$\bigoplus_{p,q} R_1(\omega_{C^*},C^*)^\Sigma_{(p+q-d+\dim C^*+1)/2}
\otimes R_1(f_C,C)_{(q-p+\dim C)/2},
$$
 should be the module over
the ring $R_0(\omega,K^*)^\Sigma/Ann([0,1])$:
$$a\cdot(b\otimes c)=\bar{a}b\otimes c,$$
for $a\in R_0(\omega,K^*)^\Sigma/Ann([0,1])$ and $(b\otimes c)$ from a component of the above
direct sum, where $\bar{a}$ is the image of $a$ induced by the projection
$R_0(\omega,K^*)^\Sigma@>>>R_0(\omega_{C^*},C^*)^\Sigma$.

We can also say about the product structure on the B-model chiral ring.
The space $R_1(f,K)\cong R_0(f,K)/Ann([0,1])$ in Conjecture~\ref{semiampleconj}, which lies
in the middle cohomology $\oplus_{p+q=d-1}H^{p,q}_{\rm st}(Y)$, should be a subring of
the  B-model chiral ring, and
$$\bigoplus_{p,q} R_1(\omega_{C^*},C^*)^\Sigma_{(p+q-d+\dim C^*+1)/2}
\otimes R_1(f_C,C)_{(q-p+\dim C)/2},
$$
 should be the module over
the ring $R_1(f,K)$, similarly to the above description in the
previous paragraph.

These ring structures are consistent  with the products on the
usual cohomology and the B-model chiral ring
$H^*(X,\bigwedge^*T_X)$ of the smooth semiample Calabi-Yau
hypersurfaces $X$ in \cite[Theorem~2.11(a,b)]{m3} and
\cite[Theorem~7.3(i,ii)]{m2}.
\end{conj}

\section{Description of string cohomology inspired by vertex algebras}
\label{section.vertex}

Here we will give yet another description of the string cohomology
spaces of Calabi-Yau hypersurfaces. It will appear as cohomology
of a certain complex, which was inspired by
the vertex algebra approach to Mirror Symmetry.

We will state the result first in the non-deformed case, and it
will be clear what needs to be done in general.
Let $K$ and $K^*$ be dual reflexive cones of dimension $d+1$ in
the lattices $M$ and $N$ respectively. We consider the subspace $\CC[L]$
of $\CC[K]\otimes \CC[K^*]$ as the span of the monomials
$[m,n]$ with $m\cdot n =0$. We also pick non-degenerate elements
of degree one $f=\sum_m f_m [m]$ and $g=\sum_n g_n [n]$ in $\CC[K]$
and $\CC[K^*]$ respectively.

Consider the space
$$
V = \Lambda^*(N_\CC)\otimes \CC[L].
$$

\begin{lem}
The space $V$ is equipped with a differential $D$ given by
$$
D:= \sum_{m} f_m \contr m \otimes(\pi_L\circ[m]) +
\sum_{n} g_n (\wedge n) \otimes(\pi_L\circ[n])
$$
where $[m]$ and $[n]$ means multiplication by the corresponding monomials
in $\CC[K]\otimes\CC[K^*]$ and $\pi_L$ denotes the natural projection
to $\CC[L]$.
\end{lem}

\begin{pf}
It is straightforward to check that $D^2=0$.
\end{pf}

\begin{thm}\label{Dcoh}
Cohomology $H$ of $V$ with respect to $D$ is naturally isomorphic
to
$$
\bigoplus_{C\subseteq K} \Lambda^{\dim C^*}C^*_\CC
\otimes R_1(f,C)\otimes R_1(g,C^*)
$$
where $C^*_\CC$ denotes the vector subspace of $N_\CC$ generated by $C^*$.
\end{thm}

\begin{pf}
First observe that $V$ contains a subspace
$$\bigoplus_{C\subseteq K} \Lambda^* N_\CC \otimes (\CC[C^\circ]\otimes
\CC[C^{*\circ}])$$
which is invariant under $D$. It is easy to calculate the cohomology
of this subspace under $D$, because the action commutes with the
decomposition $\oplus_C$. For each $C$, the cohomology of
$D$ on $\Lambda^* N_\CC \otimes (\CC[C^\circ]\otimes
\CC[C^{*\circ}])$ is naturally isomorphic to
$$\Lambda^{\dim C^*}C^*_\CC
\otimes R_0(f,C^\circ)\otimes R_0(g,C^{*\circ}),
$$
because $\Lambda^* N_\CC \otimes (\CC[C^\circ]\otimes\CC[C^{*\circ}])$
is a tensor product of the Koszul complex for $\CC[C^\circ]$ and
the dual of the Koszul complex for $\CC[C^{*\circ}]$. As a result,
we have a map
$$
\alpha:H_1\to H,~H_1:=\bigoplus_{C\subseteq K}\Lambda^{\dim C^*}C^
*_\CC
\otimes R_0(f,C^\circ)\otimes R_0(g,C^{*\circ}).
$$

Next, we observe that $V$ embeds naturally into the space
$$\bigoplus_{C\subseteq K} \Lambda^* N_\CC \otimes (\CC[C]\otimes
\CC[C^{*}])
$$
as the subspace of the elements compatible with the restriction maps.
This defines a map
$$
\beta:H\to\bigoplus_{C\subseteq K}\Lambda^{\dim C^*}C^
*_\CC
\otimes R_0(f,C)\otimes R_0(g,C^{*})=:H_2.
$$

We observe that the composition $\beta\circ\alpha$ is precisely
the map induced by embeddings $C^\circ\subseteq C$ and
$C^{*\circ}\subseteq C^*$, so its image in $H_2$ is
$$
\bigoplus_{C\subseteq K}\Lambda^{\dim C^*}C^
*_\CC \otimes R_1(f,C)\otimes R_1(g,C^{*}).
$$
As a result, what we need to show is that $\alpha$ is surjective and
$\beta$ is injective. We can not do this directly, instead, we will
use spectral sequences associated to two natural filtrations on $V$.

First, consider the filtration
$$
V=V^0\supset V^1\supset \ldots \supset V^{d+1}\supset V^{d+2}=0
$$
where $V^p$ is defined as $\Lambda^*N_\CC$ tensored with the span
of all monomials $[m,n]$ for which the smallest face of $K$ that
contains $m$ has dimension at least $p$.
It is easy to see that the spectral sequence of this filtration starts
with
$$
H_3:=\bigoplus_{C\subseteq K}\Lambda^{\dim C^*}C^
*_\CC \otimes R_0(f,C^\circ)\otimes R_0(g,C^{*}).
$$
Analogously, we have a spectral sequence from
$$
H_4:=\bigoplus_{C\subseteq K}\Lambda^{\dim C^{*}}C^
*_\CC \otimes R_0(f,C)\otimes R_0(g,C^{*\circ})
$$
to $H$, which gives us the following diagram.
$$
\begin{array}{ccccc}
   &              &   H_3      &             &     \\
   & \nearrow     &  \Downarrow& \searrow    &     \\
H_1& \rightarrow  &   H        & \rightarrow &H_2  \\
   & \searrow     & \Uparrow   & \nearrow    &     \\
   &              &   H_4      &             &
\end{array}
$$

We remark that the spectral sequences mean that $H$ is a subquotient
of both $H_3$ and $H_4$, i.e. there are subspaces $I_3^+$ and $I_3^-$
of $H_3$ such that $H\simeq I_3^+/I_3^-$, and similarly for $H_4$.
Moreover, the above diagram induces commutative  diagrams
$$
\begin{array}{ccccccccccc}
   &              &  0         &             &   &\hspace{20pt} &
   &              &  0         &             &     \\

   &              &  \downarrow&             &   &\hspace{20pt} &
   &              &  \uparrow  &             &     \\

   &              &  I_3^-     &             &   &\hspace{20pt} &
H_1& \rightarrow  &   H        & \rightarrow &H_2  \\

   &              &  \downarrow&             &   &\hspace{20pt} &
   & \searrow     &  \uparrow  & \nearrow    &     \\

   &              &  I_3^+     &             &   &\hspace{20pt} &
   &              &  I_4^+     &             &     \\

   & \nearrow     &  \downarrow& \searrow    &   &\hspace{20pt} &
   &              &  \uparrow  &             &     \\

H_1& \rightarrow  &   H        & \rightarrow &H_2&\hspace{20pt} &
   &              &  I_4^-     &             &     \\

   &              &  \downarrow&             &   &\hspace{20pt} &
   &              &  \uparrow  &             &     \\

   &              &  0         &             &   &\hspace{20pt} &
   &              &  0         &             &     \\
\end{array}
$$
with exact vertical lines. Indeed, the filtration $V^*$ induces a
filtration on the subspace of $V$
$$\bigoplus_{C\subseteq K} \Lambda^*\otimes\CC[C^\circ]\otimes
\CC[C^{*\circ}].$$
The resulting spectral sequence degenerates immediately,
and the functoriality of spectral sequences assures that there are
maps from $H_1$ as above. Similarly, the space
$$\bigoplus_{C\subseteq K} \Lambda^*\otimes\CC[C]\otimes
\CC[C^{*}]$$
has a natural filtration by the dimension of $C$ that induces
the filtration on $V$. Functoriality then gives the maps
to $H_4$.

We immediately get
$$Im(\beta) \subseteq Im(H_3\to H_2) \cap Im(H_4\to H_2)$$
which implies that
$$
Im(\beta) = Im(\beta\circ\alpha)=
\bigoplus_{C\subseteq K}\Lambda^{\dim C^*}
C^*_\CC \otimes R_1(f,C)\otimes R_1(g,C^{*}).
$$
Analogously, $Ker(\alpha)=Ker(\beta\circ\alpha)$, which shows
that
$$
\bigoplus_{C\subseteq K}\Lambda^{\dim C^*}
C^*_\CC \otimes R_1(f,C)\otimes R_1(g,C^{*})
$$
is a direct summand of $H$.

The fact that
$$Ker(\alpha)\supseteq Ker(H_1\to H_4) $$
$$=
\bigoplus_{C\subseteq K}\Lambda^{\dim C^*}C^*_\CC
\otimes Ker(R_0(f,C^\circ)\to R_0(f,C))\otimes R_0(g,C^{*\circ})
$$
implies that $I_3^-$ contains the image of this space
under $H_1\to H_3$, which is equal to
$$
\bigoplus_{C\subseteq K}\Lambda^{\dim C^*}C^*_\CC
\otimes Ker(R_0(f,C^\circ)\to R_0(f,C))\otimes R_1(g,C^{*}).
$$
Similarly, $I_3^+$ is contained in the preimage of
$$
\bigoplus_{C\subseteq K}\Lambda^{\dim C^*}C^*_\CC
\otimes R_1(f,C)\otimes R_1(g,C^{*})
$$
under $H_3\to H_2$, which is
$$
\bigoplus_{C\subseteq K}\Lambda^{\dim C^*}C^*_\CC
\otimes R_0(f,C^\circ)\otimes R_1(g,C^{*}).
$$
As a result,
$$H=\bigoplus_{C\subseteq K}\Lambda^{\dim C^*}
C^*_\CC \otimes R_1(f,C)\otimes R_1(g,C^{*}).$$
\end{pf}

\begin{rem}
If one replaces $\CC[K^*]$ by $\CC[K^*]^\Sigma$
in the definition of $\CC[L]$, then the statement and the proof
of Theorem \ref{Dcoh} remain intact. In addition, one can make
a similar statement after replacing $\Lambda^*N_\CC$ by $\Lambda^*M_\CC$
and switching contraction and exterior multiplication in the
definition of $D$. It is easy to see that the resulting complex is
basically identical, though various gradings are switched.
This should correspond to a switch between $A$ and $B$ models.
\end{rem}

We will now briefly outline the connection between Theorem \ref{Dcoh}
and the vertex algebra approach to mirror symmetry, developed
in \cite{Bvertex} and further explored in \cite{MS}.
The vertex algebra that corresponds to the N=2 superconformal
field theory is expected to be the cohomology of a lattice vertex algebra
${\rm Fock}_{M\oplus N}$, built out of $M\oplus N$, by a certain
differential $D_{f,g}$ that depends on the defining equations $f$ and
$g$ of a mirror pair.
The space $\Lambda^*(N_\CC)\otimes \CC[L]$ corresponds to a certain subspace
of ${\rm Fock}_{M\oplus N}$ such that the restriction of $D_{f,g}$ to
this subspace coincides with the differential $D$ of Theorem \ref{Dcoh}.
We can not yet show that this is precisely the chiral ring of the vertex
algebra, so the connection to vertex algebras needs to be explored further.

\section{Appendix. G-polynomials}

A finite graded partially ordered set is called Eulerian if every
its nontrivial interval contains equal numbers of elements of even and odd rank.
We often consider the poset of faces of the Gorenstein cone
$K$ over a  reflexive  polytope
$\Delta$ with respect
to inclusions. This is an Eulerian poset with
the grading given by the dimension of the face. The minimum and maximum
elements of a poset are commonly denoted by $\hat 0$ and $\hat 1$.

\begin{defn}
 \cite{stanley1} {\rm Let $P = \lbrack \hat{0}, \hat{1} \rbrack$
be an Eulerian poset of rank $d$. Define two polynomials
$G(P,t)$, $H(P,t) \in {\ZZ} [t]$ by  the following recursive rules:
$$
G(P,t) = H(P,t) = 1\;\; \mbox{\rm if $d =0$};
$$
$$
H(P,t) = \sum_{ \hat{0} <  x \leq  \hat{1}} (t-1)^{\rho(x)-1}
G(\lbrack x,\hat{1}\rbrack, t)\;\; (d>0),
$$
$$
G(P,t) =
 \tau_{ < d/2 } \left(
(1-t)H(P,t) \right) \;\;( d>0),
$$
where $\tau_{ < r }$ denotes the truncation operator
${\ZZ}\lbrack t \rbrack  \rightarrow
{\ZZ}\lbrack t \rbrack$ which is defined by
\[ \tau_{< r} \left( \sum_i a_it^i \right) = \sum_{i < r}
a_it^i. \]}
\label{Gpoly}
\end{defn}

The following lemma will be extremely useful.
\begin{lem}
For every Eulerian poset $P=[\hat 0,\hat 1]$ of positive rank there
holds
$$
\sum_{\hat 0\leq x\leq \hat 1}(-1)^{\rk[\hat 0,x]} G([\hat 0,x]^*,t)
G([x,\hat 1],t)
=\sum_{\hat 0\leq x\leq \hat 1}
G([\hat 0,x],t)G([x,\hat 1]^*,t)(-1)^{\rk[x,\hat 1]}
=0
$$
where $()^*$ denotes the dual poset. In other words, $G(\_,t)$ and $(-1)^{\rk}
G(\_^*,t)$ are inverses of each other in the algebra of functions on
the posets with the convolution product.
\label{Ginverse}
\end{lem}

{\em Proof.}
See Corollary 8.3 of \cite{stanley}.
\qed

The following polynomial invariants of Eulerian posets have been
introduced in \cite{bb}.

\begin{defn}\label{Q}
Let $P$ be an Eulerian poset of rank $d$. Define
the polynomial $B(P; u,v) \in {\ZZ}[ u,v]$
by  the following recursive rules:
$$
B(P; u,v) = 1\;\; \mbox{\rm if $d =0$},
$$
$$
\sum_{\hat{0} \leq x \leq \hat{1}}
B(\lbrack \hat{0}, x \rbrack; u,v) u^{d  - \rho(x)}
G(\lbrack x ,  \hat{1}\rbrack, u^{-1}v) = G(P ,uv).
$$
\end{defn}

\begin{lem}\label{BfromG}
Let  $P=[\hat 0,\hat 1]$ be an Eulerian poset. Then
$$
B(P;u,v) = \sum_{\hat 0\leq x \leq \hat 1} G([x,\hat 1]^*,u^{-1}v)
(-u)^{\rk \hat 1 -\rk x}
G([\hat 0,x],uv).
$$
\end{lem}

\begin{pf}
Indeed, one can sum the recursive formulas for $B([\hat 0, y])$
for all $\hat 0\leq y\leq\hat 1$ multiplied by $G([y,\hat
1]^*,u^{-1}v) (-u)^{\rk \hat 1 -\rk y}$ and use Lemma
\ref{Ginverse}.
\end{pf}

\end{document}